\documentclass[11pt]{article}
\usepackage{enumerate}
\usepackage[OT1]{fontenc}
\usepackage[usenames]{color}
\usepackage{smile}
\usepackage[colorlinks,
            linkcolor=red,
            anchorcolor=blue,
            citecolor=blue
            ]{hyperref}
\usepackage{mathrsfs}
\usepackage{fullpage}
\usepackage{hyperref}
\usepackage[protrusion=true, expansion=true]{microtype}
\usepackage{float}
\usepackage{subfigure}
\usepackage{amsfonts,amsmath,amssymb,amsthm,url,xspace}
\usepackage{tikz}
\usepackage{verbatim}
\usetikzlibrary{arrows,shapes}
\usepackage{mathtools}
\usepackage{authblk}
\usepackage{enumitem}
\usepackage{mathabx}
\usepackage{centernot}
\usepackage{comment}
\usepackage{cancel}
\usepackage[bottom]{footmisc}
\usepackage{mleftright}

\usepackage{chngcntr}
\counterwithout{equation}{section}

 \allowdisplaybreaks[1]

\newcommand{\LD}{\langle}
\newcommand{\RD}{\rangle}

\newcommand{\mQ}{\mathbb{Q}}
\newcommand{\mR}{\mathbb{R}}
\newcommand{\mE}{\mathbb{E}}

\newcommand{\mN}{\mathcal{N}}
\newcommand{\mT}{\mathcal{T}}

\newcommand{\VAR}{\text{Var}}

\newcommand{\qbeta}{q_{\hbeta}}
\newcommand{\rbeta}{r_{\hbeta}}

\newcommand{\0}{\boldsymbol{0}}
\def\ba{\boldsymbol{a}}
\def\Mid{\mathrel{\bigg|}}

\newcommand{\tte}{\tilde{e}}

\newcommand{\tbeta}{\bbeta^\star}
\newcommand{\tbetat}{{\bbeta^\star}^T}
\newcommand{\tgamma}{\bgamma^\star}
\newcommand{\bargamma}{\bar{\bgamma}}
\newcommand{\tgammat}{{\bgamma^\star}^T}

\newcommand{\ttgamma}{\tilde{\bgamma}}

\newcommand{\tG}{\bG^\star}
\newcommand{\tGt}{{\bG^\star}^T}

\newcommand{\hw}{\hat{\bw}}
\newcommand{\barw}{\bar{\bw}}
\newcommand{\hbeta}{\hat{\bbeta}}
\newcommand{\hgamma}{\hat{\bgamma}}

\newcommand{\hG}{\hat{\bG}}

\newcommand{\hf}{\hat{f}}

\newcommand{\hU}{\hat{U}}
\newcommand{\hV}{\hat{V}}

\newcommand{\tV}{{V^\star}}
\newcommand{\tU}{{U^\star}}
\newcommand{\barU}{\bar{U}}
\newcommand{\ttU}{\tilde{U}}
\newcommand{\barV}{\bar{V}}

\newcommand{\ttr}{{\tilde{r}}}

\newcommand{\ttg}{\tilde{g}}
\newcommand{\lamT}{\mT_{\lambda}}
\newcommand{\bare}{\bar{e}}
\newcommand{\hate}{\hat{e}}

\newcommand{\mEb}{\mE_{\hbeta}}
\newcommand{\mEbf}{\mE_{\hbeta, \hf}}
\newcommand{\mEf}{\mE_{\hf}}
\newcommand{\mEbx}{\mE_{\hbeta, \xb}}
\newcommand{\mEt}{\mE_{t}}

\newcommand{\IM}{\text{Im}}

\newcommand{\TR}{\text{Trace}}
\newcommand{\dist}{\text{dist}}

\def\I{\mathcal{I}}

\def\D{\mathcal{D}}

\def\wyi{\widecheck{y_i}}
\def\wy{\widecheck{y}}

\def\wSXi{\widecheck{S(\xb_i)}}

\def\wx{\widecheck{x}}
\def\wxi{\widecheck{x_i}}
\def\wzi{\widecheck{z_i}}

\def\wz{\widecheck{z}}

\def\wfxi{\widecheck{\hf(\xb_i^T\hbeta)}}

\def\ww{\widecheck{w}}
\def\wwi{\widecheck{w_i}}

\usepackage{xargs}    
\usepackage[colorinlistoftodos,prependcaption,textsize=tiny]{todonotes}
\newcommandx{\unsure}[2][1=]{\todo[linecolor=red,backgroundcolor=red!25,bordercolor=red,#1]{#2}}
\newcommandx{\change}[2][1=]{\todo[linecolor=blue,backgroundcolor=blue!25,bordercolor=blue,#1]{#2}}
\newcommandx{\info}[2][1=]{\todo[linecolor=OliveGreen,backgroundcolor=OliveGreen!25,bordercolor=OliveGreen,#1]{#2}}
\newcommandx{\improvement}[2][1=]{\todo[linecolor=Plum,backgroundcolor=Plum!25,bordercolor=Plum,#1]{#2}}

\begin{document}
\title{High-dimensional Index Volatility Models via Stein's Identity}
\author[1]{Sen Na}
\author[2]{Mladen Kolar}
\affil[1]{Department of Statistics, University of Chicago}
\affil[2]{Booth School of Business, University of Chicago}
\date{}
\maketitle

\begin{abstract}

We study the estimation of the parametric components of single and multiple index volatility models. Using the first- and second-order Stein's identities, we develop methods that are applicable for the estimation of the variance index in the high-dimensional setting requiring finite moment condition, which allows for heavy-tailed data. Our approach complements the existing literature in the low-dimensional setting, while relaxing the conditions on estimation, and provides a novel approach in the high-dimensional setting. We prove that the statistical rate of convergence of our variance index estimators consists of a parametric rate and a nonparametric rate, where the latter appears from the estimation of the mean link function. However, under standard assumptions, the parametric rate dominates the rate of convergence and our results match the minimax optimal rate for the mean index estimation. Simulation results illustrate finite sample properties of our methodology and back our theoretical conclusions.

\end{abstract}

\section{Introduction}\label{sec:1}

We consider the following index volatility model:
\begin{align}\label{mod:2}
y\mid \xb= f(\LD \tbeta, \xb\RD) + g({\tG}^T\xb)\epsilon
\end{align}
where $y$ is the response variable, $\xb\in\mR^d$ is the vector of predictors, and $\epsilon$ is a random error independent of $\xb$ with $\mE[\epsilon] = 0$ and $\mE[\epsilon^2] = 1$. In the model above, the conditional mean and variance of the response depend on the multivariate predictors only through linear projections. The unknown parts of this semi-parametric model are link functions $f:\mR\rightarrow \mR$ and $g:\mR^v\rightarrow\mR$, which are nonparametric components, and signals $\tbeta\in \mR^{d}$ and $\tG = (\tgamma_1, \ldots, \tgamma_v)\in\mR^{d\times v}$, which are parametric components satisfying $\tbetat\tbeta = 1$ and $\tGt \tG = I_v$. See \cite{Yang2017High} and \cite{Dudeja2018Learning} for a similar setup. \cite{Li1991Sliced} termed the linear space spanned by the direction of the projections as \textit{effective dimension reduction} (e.d.r.). In this paper, we focus our attention to the estimation of $\tG$.

In order to emphasize the main contribution of the work, we assume that the conditional mean of $y$ given $\xb$ follows a single index model. We note, however, that this assumption is not crucial and we will be able to estimate $\tG$ as long as the mean function can be estimated sufficiently quickly, as we illustrate later. Estimators of $f$ and $g$ only depend on the projection of predictors $\xb$ onto the e.d.r.~direction. In particular, one can apply local polynomial regression or a spline-based method on $\{y_i, {\hbeta}^T\xb_i, {\hG}^T\xb_i\}_{i=1}^n$ to estimate $f$ and $g$ once $\hbeta$ and $\hG$ are computed. Furthermore, the estimation of the nonparametric components does not depend on the ambient
dimensionality $d$ of the problem. Thus, our focus will be on
estimating parametric components in a high-dimensional setting, allowing for heavy-tailed covariates $\xb$, without using knowledge of $f$ and $g$.

Model in \eqref{mod:2} has been widely studied in the literature as it allows for flexible modeling of data without making rigid assumptions that parametric models make, while at the same time allowing for tractable estimation without suffering from the curse of dimensionality that obsesses fully nonparametric methods \citep{Robins1997curse, bach2017breaking}. When the variance function is constant and does not depend on the predictors $\xb$, the model~\eqref{mod:2} becomes the homoscedastic single index model (SIM), which plays a prominent role in econometrics and applied quantitative sciences (see, for example, \cite{Sharpe1963simplified, Collins1986Risk, Stock1988probability}). Due to its wide-ranged
applicability, a number of estimation procedures were proposed and studied (see \cite{Ichimura1993Semiparametric, Haerdle1993Optimal, Horowitz1996Direct, Xia2002adaptive, Delecroix2006semiparametric} and references therein). \cite{Li1991Sliced} developed the sliced inverse regression (SIR), which is one of the first widespread methods for estimating the e.d.r.~direction. Subsequently, a number of more advanced methods were proposed for estimating single and multiple
index models. \cite{Hristache2001Direct} estimated the
e.d.r.~direction by iteratively estimating $\tbeta$ and $f'$. \cite{Gaieffas2007Optimal} used an aggregation algorithm with local polynomial estimator to estimate $f$ at the minimax rate, while \cite{Lepski2014Adaptive} developed a procedure that adapts to the smoothness of $f$. In a setting where the dimension of the predictors, $d$, increases with the sample size, $n$, \cite{Zhu2009Nonconcave} developed a penalized inverse regression method with a nonconcave SCAD penalty and their estimator $\hbeta$ is asymptotically normal as long as $d = O(n^{1/3})$. Very recently, \cite{Dudeja2018Learning} studied the landscape of semi-parametric likelihood function of SIM by expanding $f$ in the Hermite polynomial basis. They proved that the likelihood function has no spurious local minimum. Using this result, they proposed an iterative procedure based on the gradient descent method to estimate $\tbeta$ and showed that the sample complexity of the estimator is polynomial in $d$, with the degree determined by the order of degeneracy of $f$. Besides this series of work, a number of papers have studied other index structures. For example, \cite{Carroll1997Generalized}, \cite{Mueller2001Estimation}, and \cite{Wang2010Estimation} studied partial-linear index model; \cite{Ait-Saiedi2008Cross} and \cite{Lian2011Functional} studied functional index model; and \cite{Wong2008Varying}, \cite{Xue2012Empirical}, and \cite{Ma2015Varying} studied varying-coefficient index model.

The above mentioned literature, while able to attain either $\sqrt{n}$-consistency or asymptotic normality for estimating parametric components, have two limitations. First, most of them require the predictors $\xb$ to have Gaussian or elliptically symmetric distribution. Second, they focus on estimation in a low-dimensional setting where the sample size $n$ far exceeds the dimension of predictors, $d$. More recent literature, that is closer to the approach we take in this paper, addressed these two limitations by incorporating Stein's identity into estimation of index models. \cite{Babichev2018Slice} developed a sliced inverse regression method based on Stein's identity that allows for estimation under weak conditions on the distribution of $\xb$, albeit still in a low-dimensional setting. \cite{Plan2016generalized} studied estimation of single index models in a high-dimensional setting with Gaussian design and showed that the generalized Lasso enjoys the optimal rate for the estimation of the parametric part of the model. \cite{Yang2017High} extended the above work to heavy-tailed designs, while maintaining the optimal statistical rate using the first-order Stein's identity, and further, \cite{Yang2017Estimating} developed methodology for estimation of multiple index models. \cite{Na2019High} illustrated how to estimate varying-coefficient index models. Furthermore, Stein's method has also been applied into risk estimation in Gaussian sequence model and normal approximation in recent work \citep{Chen2011Normal, Bellec2018Second}.

Allowing for conditional heteroscedasticity extends the applicability of the model even further. In financial time series, the function $g(\cdot)$ is usually interpreted as diffusion or volatility, with a long history in stochastic process, dating back to
\cite{Doob1953Stochastic}. Development of the heteroscedastic model is attributed to \cite{Engle1987Estimating}. Estimating the function $g$ does not only help in the estimation of the mean, but is interesting in its own right (see \cite{Box1974Correcting, Bickel1978Using, Box1986analysis} and references therein). \cite{Haerdle1993Optimal} first considered model~\eqref{mod:2} with $v=1$ and termed it single index volatility model, which was subsequently studied in
\cite{Xia2002Single}. \cite{Zhang2018Quasi} extended the quasi-likelihood estimator of \cite{Xia2006Asymptotic} to low-dimensional single index volatility model. \cite{Chiou2004Quasi} proposed a semiparametric quasi-likelihood approach to estimate multiple index models with purely nonparametric variance function. \cite{Klein2009semiparametric} studied a special case of a single index volatility model and built a likelihood-based estimator for unknown variance function using local smoothing. \cite{VanKeilegom2010Semiparametric} studied general
semiparametric location-dispersion models with applications to index volatility models. \cite{Fang2015Variance} proposed a two-step procedure for fitting a heteroscedastic additive partial-linear model, while \cite{Lian2015Variance} extended the method of \cite{Wang2010Estimation} for fitting a model where both mean function and variance function are in partial-linear single index form. In the literature, estimation of the variance part of index volatility models is done similarly to estimation of the mean part after properly dealing with residuals. We use these ideas to develop estimation procedures based on the Stein's identity, those allowing for estimation of volatility models in a high-dimensional setting without assuming Gaussian design.

In this paper, we address this missing study. Specifically, we
consider a generalization of the single index volatility model with $v>1$, which we call multiple index volatility model. To compare with existing literature on index volatility modeling, we focus our attention on estimation in a high-dimensional setting, which is possible under an assumption that $\tG$ is sparse. We add weak assumptions on the predictors, that allow for a wide range of heavy-tailed designs, and develop an estimator that can estimate parametric components without knowing the link function $g$, as is needed in many applications \citep{Boufounos20081, Yi2015Optimal}. In particular, we avoid iterative estimation of $\tG$ and $g$ that is common in the literature on index volatility modeling and requires some knowledge of $g$. While our estimator of $\tG$ can skip the estimation of $g$, it does rely on having a good estimator of the
conditional mean. Some necessary results concerning the mean
estimation are discussed in Section~\ref{sec:2}, while detailed
theoretical analysis is given in Appendix~\ref{sec:A}. The effect of the mean estimation will be evident from the obtained statistical convergence results, which can be decomposed into two parts: (i) nonparametric rate, originating from the estimation of $f$, and (ii) parametric rate at which we can estimate $\tG$ under the knowledge of $f$. In a high-dimensional setting, it is often the latter, parametric
part, that dominates the convergence rate, as long as $f$ is
sufficiently smooth.

The main contributions of the paper are three-fold. First, we develop a flexible method for estimating $\tG$ in the index volatility model \eqref{mod:2}, with either $v = 1$ or $v>1$. Our estimation procedure applies the first- and second-order Stein's identities on the residuals, and does not require neither sub-Gaussian or elliptical designs, nor the knowledge of the link function $g$. Second, we establish the statistical rate of convergence of the proposed estimators in a high-dimensional setting. To our knowledge, this is the first comprehensive study on high-dimensional heavy-tailed index volatility models, compared to the aforementioned literature where $\tG$ and $g$ are iteratively estimated under a low-dimensional Gaussian design. As a byproduct of our study, we also provide a result for the low-dimensional setting. Third, from technical aspect, we explicitly characterize the residuals of the local linear estimator of \eqref{mod:2}, without assuming predictors to have bounded support, which is widely used in the literature \citep[see, for example,][]{Zhu2006Empirical, VanKeilegom2010Semiparametric, Wang2010Estimation, Lian2015Variance}. This extension is critical for studying heavy-tailed designs, and bridges the gap in conditions required for estimation of $\tbeta$ and those for estimation of $f$. Specifically, while estimating $\tbeta$ can be done under weak conditions on predictors \citep{Yang2017High}, estimating $f$ requires rather stronger assumptions. In on our work, these two steps are carried jointly and, furthermore, the variance is estimated under the same, universal setup. We illustrate finite sample properties of our estimators through a series of experiments, including scenarios for which there were no suitable estimators before.

\subsection{Notation}

We summarize notations that are used throughout the paper.
We use boldface symbols to denote column vectors. For any two vectors $\ba$ and $\bb$, we use $(\ba; \bb)$ to denote a column vector obtained by stacking them together. $\be$ denotes the canonical basis of $\mR^{r}$ for some $r$ that will be clear from the context. Given an integer $k$, $[k] = \{1,2, \ldots, k\}$ denotes the index set. For any two scalars $a$ and $b$, we let $a\wedge b = \min\{a, b\}$ and $a\vee b = \max\{a, b\}$. For positive $a$ and $b$, we write $a\lesssim b$ ($a\gtrsim b$) if there exists a constant $c$ such that $a/b\leq c$ ($b/a\leq c$). We also write $a\asymp b$ if $a\lesssim b$ and $a\gtrsim b$. For a vector $\bbeta\in\mR^{d}$, we define $\|\bbeta\|_0 = |\text{supp}(\bbeta)|$. We say $\bbeta$ is $s$-sparse if $\|\bbeta\|_0\leq s$. The norm $\|\cdot\|_p$ represents either the $\ell_p$ norm of a vector or the induced $p$-norm of a matrix (for $p=2$ the norm is used without a subscript). For a matrix $A\in\mR^{m\times n}$, we let $\|A\|_*$ denote the nuclear norm, $\|A\|_F$ denote the Frobenius norm, and $\|A\|_{p,q} = \big(\sum_{j=1}^n(\sum_{i=1}^m |A_{ij}|^p)^{q/p}\big)^{1/q}$. We use $I_r$ to denote $r\times r$ identity matrix. For a random variable $v$, $\mE_{v}[\cdot] = \mE[\cdot \mid v]$ is the conditional expectation given randomness in $v$. Also, a sequence of variable $v_n$ is written as $v_n = O_P(a_n)$ if $v_n/a_n$ is stochastically bounded. Finally, $C^r(\mR)$ denotes all $r$ times continuously differentiable functions, while $\mQ^{d\times d}$ denotes all $d\times d$ orthogonal matrices.

\section{Preliminary}\label{sec:2}

In this section, we present the first- and second-order Stein's
identities that will be used as fundamental tools in our
estimation. Furthermore, we introduce the finite moment condition and some basic results on the mean estimation, under which we develop detailed estimation procedures in Sections \ref{sec:3} and \ref{sec:4}.

\subsection{Stein's Identity}

\cite{Stein1981Estimation} described the first-order Stein's identity for a Gaussian random variable, which was further extended to general random variables in \cite{Stein2004Use}. To present the first-order Stein's identity, we need the following definition.

\begin{definition}[First-order regularity condition]\label{def:1}
	
	Suppose $\Xb$ is a $\mR^d$ random vector with a differentiable density $p_{\Xb}: \mR^d\rightarrow \mR$, whose support is denoted as ${\cal X} \subseteq \mR^d$. Further, we suppose $p_{\Xb}(\bx)$ is strictly positive in the interior of $\cal X$ with $|p_{\Xb}(\bx)|\rightarrow 0$ as $\bx$ goes to the boundary. Let $S_{\Xb}: {\cal X} \rightarrow \mR^d$ be the \textit{first-order score function} defined as $S_{\Xb}(\bx) = -\nabla_{\bx}\log p_{\Xb}(\bx)$. A differentiable function $f:{\cal X} \rightarrow\mR$ together with $\Xb$ satisfies the \textit{first-order regularity condition} if both $\mE[|f(\Xb)\cdot S_{\Xb}(\Xb)|]$ and $\mE[|\nabla_{\bx} f(\Xb)|]$ exist.
	
\end{definition}

With this definition, we have the following theorem.

\begin{theorem}[First-order Stein's identity, \cite{Stein2004Use}]\label{thm:1}
	
	If function $f$ together with random vector $\Xb$ satisfies the first-order regularity condition, then we have
	\begin{align*}
	\mE[f(\Xb)\cdot S_{\Xb}(\Xb)] = \mE[\nabla_{\bx} f(\Xb)].
	\end{align*}
	
\end{theorem}

In order to generalize to the second-order identity, we define the second-order regularity condition.

\begin{definition}[Second-order regularity condition]
	
	Suppose the same conditions as in Definition \ref{def:1} hold. Let $H_{\Xb}: {\cal X}\rightarrow \mR^{d\times d}$ be the \textit{second-order score function} defined as $H_{\Xb}(\bx) = \nabla_{\bx}^2 p_{\Xb}(\bx)/p_{\Xb}(\bx)$. A twice differentiable function $f: {\cal X}\rightarrow \mR$ together with $\Xb$ satisfies the \textit{second-order regularity condition} if both $\mE[|f(\Xb)\cdot H_{\Xb}(\Xb)|]$ and $\mE[|\nabla_{\bx}^2 f(\Xb)|]$ exist.
	
\end{definition}

\begin{theorem}[Second-order Stein's identity, \cite{Janzamin2014Score}]\label{thm:2}
	
	If function $f$ together with random vector $\Xb$ satisfies the second-order regularity condition, then we have
	\begin{align*}
	\mE[f(\Xb)\cdot H_{\Xb}(\Xb)] = \mE[\nabla_{\bx}^2 f(\Xb)].
	\end{align*}
	
\end{theorem}

In what follows, we will omit the subscript in $\nabla_{\bx}f, S_{\Xb}, H_{\Xb}$ and write $\nabla f, S, H$ whenever it's clear from the context. It is easy to see that when $\Xb\sim \mN(\0, I_{d})$, then $S(\Xb) = \Xb$, $H(\Xb) = \Xb\Xb^T - I_d$. Furthermore, by the above two theorems, we get
\begin{align*}
\mE[f(\Xb) \cdot \Xb] = \mE[\nabla f(\Xb)] \text{\ \  and\ \ } \mE[f(\Xb)\cdot (\Xb\Xb^T - I_d)] = \mE[\nabla^2 f(\Xb)],
\end{align*}
if $f$ satisfies both regularity conditions. The regularity conditions are fairly mild and~are required in the literature on Stein-based estimators. See, for example, \cite{Yang2017High, Babichev2018Slice, Na2019High} and references~therein. In addition to the regularity conditions above, we will need a moment assumption.

\begin{assumption}[Finite moment assumption]\label{ass:1}
	
	We say finite $p$-th moment assumption holds for the model \eqref{mod:2}, if $\mE[|\epsilon|^p]<\infty$ and there exists constant $M_p>0$ such that $\forall j,k\in[d]$ and unit vector $\vb\in\mR^d$,
	\begin{align*}
	\mE[|\xb^T\vb|^p] \vee \mE[|f(\tbetat\xb)|^p] \vee \mE[|g(\tGt\xb)|^p] \vee\mE[|S(\xb)_j|^p]\vee \mE[|H(\xb)_{jk}|^p]\leq M_p.
	\end{align*}
	Furthermore, we have
	\begin{align*}
	\mE[|y|^p]\lesssim \mE[|f(\LD \xb, \tbeta\RD)|^p] + \mE[|g(\tGt\xb)|^p]\cdot\mE[|\epsilon|^p]\lesssim M_p.
	\end{align*}
	
\end{assumption}

In the above assumption, we assume that $\mE[|\epsilon|^p]$ is a constant and do not keep track of it. On the other hand, we explicitly keep track of the quantity $M_p$. Although the above assumption does not explicitly put restrictions on the tails of $\xb$, it does allow for certain types of heavy-tailed designs, including Gamma and $t$-distribution. Furthermore, when the predictor $\xb$ has $i.i.d$ entries, $\mE[|H(\xb)_{jk}|^p]$ is bounded as long as the $p$-th moment of $S(\xb)_j$ and its derivative are bounded.

\subsection{Mean Estimation}

Our estimator for $\tG$ in the model \eqref{mod:2} relies on a good estimator of the conditional mean. Since the variance estimation procedure does not depend on the specific form of the conditional mean function, in order to simplify the presentation, we first consider the following model,
\begin{align}\label{mod:1}
y \mid \xb= f(\xb) + g(\tGt\xb)\epsilon,
\end{align}
where $\xb\in\mR^d$ is the predictor vector, $\epsilon$ is noise with $\mE[\epsilon | \xb] = 0$, and $f(\cdot)$ is an unknown function that is not necessarily of the index form. While this model is not suitable in a high-dimensional setting, it helps us illustrate the main requirements on the conditional mean estimator. Detailed estimation procedure, assumptions, and convergence results for the index model \eqref{mod:2} are provided in Appendix \ref{sec:A}.

Under the model \eqref{mod:1} with $\xb$ belonging to a compact set, a number of standard nonparametric methods can be used to estimate $f$, such as local polynomial regression and smoothing spline. Suppose we use $n$ independent samples, say $\D = \{y_i, \xb_i\}_{i=1}^n$ to
estimate $\hf$, under suitable regularity conditions (see, for
example, the Condition 1 in \cite{Fan1993Local} for one-dimensional case), the pointwise mean squared error can be upper-bounded as
\begin{align}\label{cond:1}
\mE[|\hf(\xb) - f(\xb)|^2 \mid \xb = \xb_0] \leq  e_f(\xb_0, n, d)
\end{align}
for some error function $e_f(\xb_0, n, d)$ depending on the evaluation point $\xb_0$, dimension $d$ and sample size $n$. In particular, when $f\in \bSigma(k, L)$ where $\bSigma(k, L)$ denotes the H\"older class indexed by $k$ and $L$ (see Definition 1.2 in \cite{Tsybakov2009Introduction}), the integrated mean squared error satisfies
\begin{align}\label{cond:2}
\mE[|\hf(\xb) - f(\xb)|^2] \leq \mE[e_f(\xb, n, d)] \leq \Upsilon\cdot n^{-\frac{2k}{2k+d}},
\end{align}
for some constant $\Upsilon$, which is also the minimax rate
\citep{Gyoerfi2002distribution}.

Different from the above discussed mean estimation, in order to have precise variance information for a given $\hf$, we require a slightly stronger result on $\hf$. Suppose $W(\xb)$ is an entry of either the first- or second-order score variable. We require that the weighted mean squared error, for the given $\hf$, is well controlled. In particular, $\forall 0<\delta<1$,
\begin{align}\label{equ:Cauchy}
\mE\big[|\hf(\xb) - f(\xb)|^2\cdot W(\xb) \mid \D\big]\leq &\sqrt{\mE[|\hf(\xb) - f(\xb)|^4 \mid \D]}\cdot \underbrace{\sqrt{\mE[W(\xb)^2]}}_{\text{bounded}} \nonumber\\
\lesssim &\sqrt{\mE[|\hf(\xb) - f(\xb)|^4]}/\delta\eqqcolon \tte_f(n, d)/\delta,
\end{align}
where the second inequality is due to the Markov's inequality and
holds with probability $1-\delta$. To have the weighted mean squared error bounded, we require the mean estimator to satisfy
\begin{align}\label{cond:3}
\tte_f(n, d) = \sqrt{\mE[\bare_f(\xb, n, d)]}<\infty,
\end{align}
where
\begin{align}\label{cond:1new}
\bare_f(\xb_0, n, d) \coloneqq \mE[|\hf(\xb) - f(\xb)|^4 \mid \xb = \xb_0].
\end{align}
Compared to the bounded second moment in \eqref{cond:1}, here we require the fourth moment to be bounded in \eqref{cond:1new}, due to the application of the Cauchy-Schwarz in \eqref{equ:Cauchy}. When the covariate vector $\xb$ is supported on a compact set and its density
is bounded away from zero, one can simply show that $e_f(\xb, n, d)$ is uniformly upper bounded and, in fact, converges to zero at the rate of $n^{-\frac{2k}{2k+d}}$, and as a result has all the moments bounded.

In order to emphasize the main contribution, which is the variance estimation, we use the local linear estimator proposed in \cite{Fan1993Local} and derive an explicit formula for
$\bare_f(\xb, n, d)$ under the model~\eqref{mod:2} in
Appendix~\ref{sec:A}. We prove that \eqref{cond:3} holds under a tail condition on $\bare_f(\xb, n, d)$, which is satisfied for any compact designs, as well as for any link functions $f$ with appropriate decay properties. We note that an alternative proof technique is possible using uniform convergence result of $\hf$, see \cite{Hansen2008Uniform}, which, however, would require different regularity conditions.

Under \eqref{cond:3}, the following theorem provides a result on the mean estimation that will allow us to estimate $\tG$ in model \eqref{mod:1}.

\begin{theorem}\label{thm:mean}
	
	Suppose there is an estimator $\hat f$ of the conditional mean under the model \eqref{mod:1} which is calculated from $n$ samples $\D = \{y_i, \xb_i\}_{i\in[n]}$ and satisfies condition \eqref{cond:3}. Let $W(\xb)$ be either $S(\xb)$ or $H(\xb)$ and assume that each entry of $W(\xb)$ has a finite $2$-nd moment bounded by $M_2$. Then, for any $0<\delta<1$, we have
	\begin{align}\label{cond:4}
	P\bigg(\biggl\|\mE\big[|\hf(\xb) - f(\xb)|^2\cdot W(\xb) \mid \D\big]\biggr\|_{\infty}\geq \frac{\sqrt{M_2}\cdot\tte_f(n,d)}{\delta} \bigg)\leq \delta.
	\end{align}
	
\end{theorem}

The following result is an immediate corollary for the index model \eqref{mod:2}.

\begin{corollary}\label{cor:1}
	
	Suppose $\hbeta$ and $\hf$ are estimators of $\tbeta$ and $f$ under the model \eqref{mod:2}, calculated from two independent sample sets $\D_1$ and $\D_2$ with size $n$ for each. We define the mean quartic error as
	\begin{align}\label{cond:5}
	\bare_f(\hbeta^T\xb, n, 1)\coloneqq & \mE\big[|\hf(\LD\hbeta, \xb\RD) - f(\LD\tbeta, \xb\RD)|^4 \mid \xb, \D_1\big]
	\end{align}
	and assume
	\begin{equation}\label{cond:6}
	\begin{aligned}
	&\hate_f(\hbeta, n, 1) \coloneqq \mE\big[\bare_f(\hbeta^T\xb, n, 1) \mid \D_1\big]<\infty,\\
	&P\big(\sqrt{\hate_f(\hbeta, n, 1)} \geq \tte_{f, \delta}(n, 1) \big)\leq \delta, \text{\ \ \ } \forall 0<\delta<1,
	\end{aligned}
	\end{equation}
	for some rate $\tte_{f, \delta}(n, 1)$. Let $W(\xb)$ be either $S(\xb)$ or $H(\xb)$ and assume that each entry of $W(\xb)$ has a finite $2$-nd moment bounded by $M_2$. Then we have
	\begin{align}\label{cond:7}
	P\bigg(\biggl\|\mE\big[|\hf(\LD\hbeta, \xb\RD) - f(\LD\tbeta, \xb\RD)|^2\cdot W(\xb) \mid \D_1, \D_2\big]\biggr\|_{\infty}\geq \frac{\sqrt{M_2}\cdot\tte_{f, \delta}(n,1)}{\delta} \bigg)\leq 2\delta,
	\end{align}
	where probability is taken over randomness in $\D_1$ and $\D_2$.
	
\end{corollary}

Estimation procedures we develop in Section \ref{sec:3} and
\ref{sec:4} for $\tG$ assume that the mean estimation satisfies
\eqref{cond:7}. In Appendix \ref{sec:A}, we will provide a simple
estimator that indeed satisfies \eqref{cond:7} for completeness. In particular, we show that $\hbeta$ can be estimated using the approach proposed in \cite{Yang2017High}, while $\hf$ can be estimated by local linear regression \citep{Fan1993Local}. However, note that a number of alternative procedures, such as smoothing splines \citep{Boor2001practical, Green1993Nonparametric}, wavelets \citep{Johnstone2011Gaussian, Mallat2009wavelet}) could be used, since under standard assumptions the quantity in \eqref{cond:5} can be uniformly bounded over evaluation points. We show that the condition \eqref{cond:6} follows from an explicit formula for $\bare_f(\hbeta^T\xb, n, 1)$. Moreover, when $f\in\bSigma(2, L)$, Theorem \ref{Supthm:2} shows that
$\tte_{f, \delta}(n, 1) \asymp n^{-4/5}$. Our analysis recovers the
existing results on estimating $f$ under model \eqref{mod:2} when $\xb$ is in a compact set, however, a more careful analysis is needed when $\xb$ is heavy-tailed.

In the following two sections, we assume existence of the estimators of $\hbeta$ and $\hf$ under model~\eqref{mod:2}, which satisfy the error rate in \eqref{cond:7}. Furthermore, to simplify the presentation of the paper, we assume that estimation of $\tG$ is done on an independent sample set with size $n$, which ensures the independence of $\hG$ from $\hbeta$ and $\hf$. This can be achieved through data splitting and will not affect the statistical rate of convergence, but only the constants.

\section{Single Index Volatility Model}\label{sec:3}

We start our analysis by focusing on single index volatility models, which are a sub-class of model in \eqref{mod:2} with $v = 1$. In particular, we focus on the following model
\begin{align}\label{mod:4}
y\mid \xb = f(\LD \xb, \tbeta\RD) + g(\LD \xb, \tgamma\RD)\epsilon,
\end{align}
and develop a procedure for estimating $\tgamma$. As discussed in Section \ref{sec:2}, we assume existence of estimators $\hbeta$ and $\hat f$ that satisfy \eqref{cond:7}. We present our estimators based on the first- and second-order Stein's identities in the following two subsections.

\subsection{First-order estimation}

Suppose the function $g^2(\LD \xb, \tgamma\RD)$ together with $\xb$ satisfies the first-order regularity condition, then
\begin{align}\label{equ:1}
\mE\big[\big(y - f(\LD\xb, \tbeta\RD)\big)^2S(\xb)\big] = \mE [\epsilon^2g^2(\LD \xb, \tgamma\RD)S(\xb)] = 2\mu_1\tgamma \eqqcolon\ttgamma,
\end{align}
where $\mu_1 = \mE[g(\LD\xb, \tgamma\RD)g'(\LD\xb, \tgamma\RD)]$. Note that whenever  $\mu_1\neq 0$, the line spanned by $\tgamma$ is identifiable from $\ttgamma$. In particular, one can estimate $\pm\tgamma$ by normalizing the estimator of $\ttgamma$. If we further assume that $\mu_1>0$ or the first entry of $\tgamma$ is positive, then one can fully identify $\tgamma$ from $\ttgamma$ \citep{Xia2006Asymptotic, Wang2010Estimation}. We take a different approach and, in order to avoid issues with normalization, use the following distance
\begin{align}\label{equ:cos}
{\rm dist}(\hgamma, \tgamma) = 1 - \frac{|\LD \hgamma, \tgamma\RD|}{\|\hgamma\|_2\|\tgamma\|_2},
\end{align}
as a surrogate of $\|\hgamma - \ttgamma\|_2^2$, to quantify the convergence rate for the first-order estimator. We will estimate $\ttgamma$ by replacing the left hand side in \eqref{equ:1} with its truncated empirical counterpart. We first introduce the truncation notation.

\begin{definition}[Truncation function]
	
	For a scalar $v \in \mR$, the truncation function is defined as $\Psi_{\tau}(v) = v\cdot\pmb{1}_{\{|v|\leq \tau\}}$. For a vector or matrix $\bv$, the truncation function $\Psi_{\tau}(\bv)$ is applied elementwise.
	
\end{definition}

In a low-dimensional setting, we estimate $\ttgamma$ as
\begin{align}\label{est:1}
\hgamma_1 = \frac{1}{n}\sum_{i=1}^n\rbr{\Psi_{\tau}(y_i) - \Psi_{\tau}(\hat f(\xb_i^T\hbeta))}^2\cdot \Psi_{\tau}(S(\xb_i)).
\end{align}
The following theorem gives us its statistical convergence rate.

\begin{theorem}[Low-dimensional first-order estimator]\label{Sthm:1}
	
	Suppose the function $g^2(\LD \xb, \tgamma\RD)$ together with $\xb$ satisfies the first-order regularity condition. Furthermore, suppose Assumption \ref{ass:1} holds with $p\geq 6$, $f$ is Lipschitz continuous and $\mu_1\neq 0$. Then, for any $0<\delta<1$, there exist constants $N_{\delta}$ (depending on $\delta$) and $\Upsilon$ such that the estimator \eqref{est:1} with $\tau = \Upsilon\rbr{\frac{nM_6}{\log(d/\delta)}}^{1/6}$ satisfies
	\begin{align*}
	P\bigg(\|\hgamma_1 - \ttgamma\|_{\infty}\leq \Upsilon\big(\sqrt{\frac{M_6\log(d/\delta)}{n}} + \frac{\sqrt{M_6}\cdot\tte_{f, \delta}(n,1)}{\delta} \big) \bigg)\geq 1-\delta,
	\end{align*}
	for all $n\geq N_{\delta}$. In addition, by the fact that $\|\hgamma_1 - \ttgamma\|_2\leq \sqrt{d}\|\hgamma_1 - \ttgamma\|_{\infty}$, if the mean estimation in Theorem \ref{Supthm:1} and \ref{Supthm:2} are satisfied, then we have
	\begin{align*}
	\dist(\hgamma_1, \tgamma)= O_P\rbr{\frac{d\log d}{\mu_1^2n}}.
	\end{align*}
	
\end{theorem}

Lipschitz continuity on the conditional mean function, also assumed in \cite{Xia2006Asymptotic} and \cite{Wang2010Estimation}, is only required to hold in a small neighborhood of $\tbeta$. See Assumption C2 in \cite{Xia2006Asymptotic}, for example. It guarantees that $f(\xb^T\hbeta)$ and $f(\xb^T\tbeta)$ are close to each other for a given estimator $\hbeta$ and, furthermore, $\mEbf[|\hf(\xb^T\hbeta)|^6]\lesssim M_6$ with high
probability. Lipschitz condition can be imposed on $f'$ or
$f''$ if $\xb$ has higher order moments \citep{Zhang2018Quasi}.

The convergence rate consists of two parts: parametric rate and
nonparametric rate. When $f\in \bSigma(2, L)$, Theorem \ref{Supthm:2} shows that $\tte_{f, \delta}(n,1)\asymp n^{-4/5}$ and therefore the parametric rate is the dominant one. Generally, when a one-dimensional function $f\in\bSigma(k, L)$, we have that $\tte_{f, \delta}(n,1)\asymp n^{-\frac{2k}{2k+1}}$, and the dominant term will be the parametric rate. A similar rate has also been established in \cite{Dudeja2018Learning}, where authors analyzed the landscape of nonconvex, semiparametric likelihood function. However, their analysis requires Gaussian design. It is not clear if their landscape analysis can be extended to general distributions.

In a high-dimensional setting, estimation of $\tgamma$ is possible under additional structural assumptions on the unknown vector. It is common to assume that $\tgamma$ is sparse and satisfies $\|\tgamma\|_0\leq s$. Under this assumption, we propose the following $\ell_1$-penalized estimator
\begin{align}\label{est:2}
\hgamma_2 = \arg\min_{\bgamma}\ \frac{1}{2}\|\bgamma\|^2 - \LD \bgamma, \hgamma_1 \RD + \lambda\|\bgamma\|_1.
\end{align}
It is well known that $\hgamma_2$ can be obtained by soft-thresholding
$\hgamma_1$ as $\hgamma_2 = \phi_{\lambda}(\hgamma_1)$, where the soft-thresholding function $\phi_{\lambda}(v) = (1- \lambda/|v|)_+\cdot v$  is applied elementwise. We have the following convergence result.

\begin{theorem}[High-dimensional first-order estimator]\label{Sthm:2}
	
	Suppose the conditions of Theorem \ref{Sthm:1} are satisfied and further suppose $\|\tgamma\|_0\leq s$. Then, for the same constants $\Upsilon$ and $N_{\delta}$, the estimator \eqref{est:2} with $\hgamma_1$ as in Theorem \ref{Sthm:1} and
	\begin{align*}
	\lambda\geq 2\Upsilon\bigg(\sqrt{\frac{M_6\log(d/\delta)}{n}} + \frac{\sqrt{M_6}\cdot \tte_{f, \delta}(n,1)}{\delta}\bigg)
	\end{align*}
	satisfies
	\begin{align*}
	P\rbr{\|\hgamma_2 - \ttgamma\|_2\leq 3\sqrt{s}\lambda\text{\ \ and\ \ } \|\hgamma_2 - \ttgamma\|_1\leq 12s\lambda}\geq 1-\delta,
	\end{align*}
	for all $n\geq N_{\delta}$. In addition, if the mean estimation in Theorem \ref{Supthm:1} and \ref{Supthm:2} are satisfied, we have
	\begin{align*}
	\dist(\hgamma_2, \tgamma)= O_P\rbr{\frac{s\log d}{\mu_1^2n}}.
	\end{align*}
	
\end{theorem}

The above two theorems show that $\tgamma$ can be estimated at a parametric rate in a low-dimensional setting and at the rate
$\sqrt{s\log d/n}$ in a high-dimensional setting. These rates are
minimax optimal when estimating mean signal $\tbeta$ in homoscedastic index model \citep{Lin2017optimality}. The results only hold asymptotically due to the estimation of the link function $f$.

\subsection{Second-order estimation}
\label{sec:sivm:sec-order-estimation}

In this section, we develop the second-order estimation procedure for $\tgamma$. Though the first-order estimator is easy to compute and has good statistical convergence rate, it has been observed in the literature that second-order estimators are more robust and the regularity condition allows for estimation under a wider class of functions \citep{Babichev2018Slice}.

Suppose the function $g^2(\LD \xb, \tgamma\RD)$ together with $\xb$ satisfies the second-order regularity condition. Under the model \eqref{mod:4}, we have
\begin{align}\label{equ:3}
\tU\coloneq \mE[(y - f(\LD\xb, \tbeta\RD))^2H(\xb)] = \mE [\epsilon^2g^2(\LD \xb, \tgamma\RD)H(\xb)]= 2\mu_2\tgamma\tgammat,
\end{align}
where
$\mu_2 = \mE[(g'(\xb^T\tgamma))^2]+
\mE[g(\xb^T\tgamma)g''(\xb^T\tgamma)]$. Suppose $\mu_2\neq 0$, one strategy for estimating $\pm\tgamma$ is based on estimating the matrix $\tU$ and then extracting its leading eigenvector. In a low-dimensional setting, this strategy leads to our second-order estimator, which is defined as
\begin{align}\label{est:4}
\hgamma_3 \in \arg\max_{\|\bgamma\|_2\leq 1} \bigg|\bgamma^T\bigg(\underbrace{\frac{1}{n}\sum_{i=1}^n\big(\Psi_{\tau}(y_i) - \Psi_{\tau}(\hf(\xb^T\hbeta))\big)^2\cdot\Psi_{\tau}(H(\xb_i))}_{\hU}\bigg)\bgamma\bigg|.
\end{align}
The following theorem establishes its rate of convergence.

\begin{theorem}[Low-dimensional second-order estimator]\label{Sthm:4}
	
	Suppose the function $g^2(\LD \xb, \tgamma\RD)$ together with $\xb$ satisfies the second-order regularity condition.  Furthermore, suppose Assumption~\ref{ass:1} holds with $p\geq 6$, $f$ is Lipschitz continuous and $\mu_2 \neq 0$. Then, for any $0<\delta<1$, there exist constants $N_{\delta}$ (depending on $\delta$) and $\Upsilon$ such that the estimator \eqref{est:4} with $\tau = \Upsilon\rbr{\frac{nM_6}{\log(d/\delta)}}^{\frac{1}{6}}$ satisfies
	\begin{align*}
	P\rbr{\min_{\iota = \pm 1}\|\iota \hgamma_3 - \tgamma\|_2\leq  \frac{\Upsilon}{\mu_2}\rbr{d\sqrt{\frac{M_6\log(d/\delta)}{n}} + \frac{d\sqrt{M_6}\cdot\tte_{f, \delta}(n,1)}{\delta} } }\geq 1-\delta,
	\end{align*}
	for all $n\geq N_{\delta}$. In addition, if the mean estimation in
	Theorem \ref{Supthm:1} and \ref{Supthm:2} are satisfied, then we have
	\begin{align*}
	\min_{\iota = \pm 1}\|\iota \hgamma_3 - \tgamma\|_2 = O_P\rbr{\frac{d}{\mu_2}\sqrt{\frac{\log d}{n}}}.
	\end{align*}
	
\end{theorem}

Theorem \ref{Sthm:4} suggests that the estimator $\hgamma_3$ attains $\sqrt{n}$-consistency in a low dimensional setting, which matches the optimal parametric rate. Our rate is the same as the one established in \cite{Babichev2018Slice}, while we relax their assumptions from sub-Gaussian condition to finite moment condition. \cite{Lin2018consistency} improved the dependence on $d$ in the rate by applying sliced inverse regression and letting sizes of slices tend to infinity. It is unclear how to improve the rate for heavy-tailed designs or without doing slicing. In the following remark, we discuss an alternative estimator in the Gaussian setup with improved rate. The modified estimator uses a soft-truncated estimator of $\tU$, as the hard-truncated estimator, $\hU$, is not suitable due to its loose error bound in $\|\cdot\|_2$.

	\begin{remark}
		
		For the model \eqref{mod:4}, suppose $\xb$ is Gaussian, $y$ and $f(\xb^T\tbeta)$ have bounded $6$-th moment. Define the following soft-truncated estimator
		\begin{align*}
		\hU' = \frac{1}{n\kappa}\sum_{i = 1}^n \Phi\big(\kappa\cdot(y_i - \hf(\xb^T_i\hbeta))^2H(\xb_i)\big),
		\end{align*}
		where $\Phi(\cdot)$ is defined in \cite{Minsker2018Sub} (equation (3.5)) and $\kappa$ is a truncation threshold.  Similar to the proof of Lemma 10 in \cite{Na2019High}, we can show that $\|\hU' - \tU\|_2\lesssim \sqrt{d\log d/n}$ with
		$\kappa\asymp \sqrt{\log d/(nd)}$. Hence, extracting the leading eigenvector based on $\hU'$ will result in the rate $\sqrt{d\log d/n}$ for estimating $\tgamma$. Comparing with the second-order estimator in \cite{Dudeja2018Learning}, we do not require the error term to be Gaussian, nor $f$ to be bounded. In addition, our estimator has a closed form, while \cite{Dudeja2018Learning} relies on an iterative algorithm. Finally, we note that although $\hU'$ attains a better dependence on $d$, it is not suitable in high dimensions, in contrast to $\hU$, due to its loose error bound in $\|\cdot\|_{\infty, \infty}$.
		
	\end{remark}

Based on $\hU$ defined in \eqref{est:4}, we proceed to estimating $\tgamma$ in high dimensions. Our estimator is built on the optimization algorithm that was proposed as a convex relaxation for sparse PCA problem \citep{Vu2013Fantope}. Given a symmetric matrix $A$, tuning parameter $\lambda$, and an integer $r$, we denote $\mT_{\lambda}(A, r)$ to be the optimal solution of the following
optimization program
\begin{equation}\label{equ:4}
\begin{aligned}
\lamT(A, r) = \arg\max_{V}\ \  &\LD V, A\RD - \lambda\|V\|_{1,1},\\
\text{s.t. } &0\preceq V\preceq I_d,\ \ \TR(V) = r.
\end{aligned}
\end{equation}
The constraint set in \eqref{equ:4} is called the Fantope of order
$r$, which is the convex hull of rank-$r$ projection matrices
\citep{Vu2013Fantope}. The tuning parameter $r$ controls the number of eigenvectors we aim to estimate, while $\lambda$ controls the overall sparsity of eigenvectors. Let $\hV = \lamT(\hU, 1)$ where $\hU$ is defined in (\ref{est:4}). Our high-dimensional second-order estimator is further defined as
\begin{align}\label{est:5}
\hgamma_4 \in \arg\max_{\|\bgamma\|_2\leq 1} |\bgamma^T\hV\bgamma|.
\end{align}

\begin{theorem}[High-dimensional second-order estimator]\label{Sthm:5}
	
	Suppose the conditions of Theorem \ref{Sthm:4} are satisfied and further suppose $\|\tgamma\|_0\leq s$. Then, for the same constants $\Upsilon$ and $N_\delta$, the estimator \eqref{est:5} with $\hU$ as in Theorem \ref{Sthm:4} and
	\begin{align*}
	\lambda\geq \Upsilon\bigg(\sqrt{\frac{M_6\log(d/\delta)}{n}} + \frac{\sqrt{M_6}\cdot\tte_{f, \delta}(n,1)}{\delta}\bigg)
	\end{align*}
	satisfies
	\begin{align*}
	P\bigg(\min_{\iota = \pm 1}\|\iota \hgamma_4 - \tgamma\|_2\leq \frac{4\sqrt{2}s\lambda}{\mu_2}\bigg)\geq 1-\delta,
	\end{align*}
	for all $n\geq N_{\delta}$. In addition, if the mean estimation in Theorem \ref{Supthm:1} and \ref{Supthm:2} are satisfied, we have
	\begin{align*}
	\min_{\iota = \pm 1}\|\iota \hgamma_4 - \tgamma\|_2 =O_P\rbr{\frac{s}{\mu_2}\sqrt{\frac{\log d}{n}}}.
	\end{align*}
	
\end{theorem}

From Theorem~\ref{Sthm:5}, we see that the second-order estimator $\hgamma_4$ achieves $s\sqrt{\log d/n}$ rate of convergence in a high-dimensional setting. Compared to the first-order estimators, the high-dimensional rate has an extra $\sqrt{s}$ factor, which comes from the convex relaxation programming we adopted. Our rate matches the one in \cite{Vu2013Fantope}, even though the estimation is done on
truncated data due to heavy-tailedness. As discussed in \cite{Vu2013Fantope}, this $\sqrt{s}$ factor is necessary for having a polynomial time procedure when doing sparse PCA under some cases. On the other hand, the success of existing sparse PCA algorithms relies on a sharp concentration of the restricted operator norm defined as
	\begin{align*}
	\|\hU - \tU\|_{2,s} = \sup_{\|\bv\|_0\leq s, \|\bv\|_2=1}|\bv^T(\hU - \tU)\bv |.
	\end{align*}
	To the best of our knowledge, when each entry of $\hU$ only has finite moment, it is an open question whether the restricted operator norm above has the rate $\sqrt{s\log d/n}$. Under a Gaussian (or sub-Gaussian) setup, \cite{Vu2013Minimax} proved that the restricted operator norm concentrates with rate $\sqrt{s\log d/n}$ and, hence, one can apply, for example, the two-stage algorithm in \cite{Wang2014Tighten} to attain the optimal $\sqrt{s\log d/n}$ rate. Despite the worse rate, the second-order estimator requires $\mu_2\neq 0$, which is a milder condition than $\mu_1\neq 0$.  For example, if $\xb^T\tbeta$ has a symmetric distribution and $g(x) = x^k$ for some $k$, then $\mu_1 = 0$, while $\mu_2\neq0$. Therefore, each estimator has its own advantages.

So far, we have investigated first- and second-order estimators under model \eqref{mod:4} and derived asymptotic results on their convergence in different settings. In the following, we will discuss an estimation procedure that can obtain a finite sample result in a setting where the mean and variance index are approximately orthogonal.

\subsection{Estimation under orthogonality}

In the previous two subsections, we discussed estimators of $\tgamma$ that require estimation of $\tbeta$ and~$f$. It is interesting to point out that estimation of $\tgamma$ is possible without estimating $f$ in a certain setting. To illustrate this, we consider the model \eqref{mod:4} in a high-dimensional setting. We have following two observations. First, if $\tbeta$ and $\tgamma$ are suitably orthogonal, then our estimation procedure should take advantage of this property. Second, two sparse vectors in high dimensions are
orthogonal to each other with high probability.

We illustrate the second point from a Bayesian point of view. Suppose that $\tbeta$ and $\tgamma$ are drawn from a prior that puts a lot of mass on $s$-sparse vectors. For example, consider the following mixture distribution from which each entry of $\tbeta$ and $\tgamma$ are drawn independently
\begin{align*}
\beta_i^\star, \gamma_i^\star\sim (1 - \pi)\cdot\delta_0 + \pi\cdot\mN(0,1), \text{\ \ \ } \forall i\in[d]
\end{align*}
where $\delta_0$ is the Dirac function, putting all mass on $0$, and $\pi= s/d$. Such a mixture distribution has been widely used in high-dimensional sparse parameter estimation
\citep{Johnstone2004Needlesa}, variable selection \citep{Mitchell1988Bayesian, Ishwaran2005Spike}, multi-task learning \citep{Titsias2011Spike}, and Bayesian multiple testing
\citep{Scott2006exploration}. Under this prior, we have that
$\tbetat\tgamma = 0$ with high probability, since
\begin{align*}
P(\tbetat\tgamma = 0)\geq\prod_{i\in[d]}P(\beta_i^\star\gamma_i^\star = 0)= (1 - \pi^2)^d \geq 1 - \frac{s^2}{d}.
\end{align*}
Next, we show that when $\tbetat\tgamma = 0$, we can estimate $\tgamma$ without estimating $f$.

We start with an estimator based on the first-order Stein's
identity. Suppose $f^2(\LD \xb, \tbeta\RD)$ and $g^2(\LD\xb, \tgamma\RD)$ together with $\xb$ satisfy the first-order
regularity condition, then
\begin{align}\label{equ:2}
\mE[y^2S(\xb)] = \mE [f^2(\LD \xb, \tbeta\RD)S(\xb)] + \mE [g^2(\LD \xb, \tgamma\RD)S(\xb)] = 2\eta_1\tbeta + \ttgamma,
\end{align}
where $\eta_1 = \mE[f( \xb^T\tbeta)f'(\xb^T\tbeta)]$. We utilize
\eqref{equ:2} to obtain our estimator. First, we can use the procedure of \cite{Yang2017High} to estimate $\tbeta$. Note that other estimators are also possible, as long as $\hbeta$ satisfies the convergence rate in Theorem~\ref{Supthm:1}. Next, given user specified thresholds $\tau$ and $\kappa$, we define
$\barw =
\frac{1}{n}\sum_{i=1}^{n}\Psi_{\tau}(y_i)^2\cdot\Psi_{\tau}(S(\xb_i))$
and its soft-thresholded version $\hw = \phi_{\kappa}(\barw)$. Our estimator is then given as
\begin{align}\label{est:3}
\hgamma_5 =  \phi_{\lambda}(\barw - \LD \hbeta, \hw\RD\cdot\hbeta),
\end{align}
where $\lambda$ is a user specified parameter that controls the
sparsity of the estimator. Its statistical convergence rate is
provided in next theorem.

\begin{theorem}[First-order orthogonal estimation]\label{Sthm:3}
	
	Suppose $f^2(\LD \xb, \tbeta\RD)$ and $g^2(\LD\xb, \tgamma\RD)$ together with $\xb$ satisfy the first-order regularity condition, $\|\tbeta\|_0 \vee \|\tgamma\|_0\leq s$, $\LD \tbeta,\tgamma\RD=0$, and $\mu_1\neq 0$ (defined in \eqref{equ:1}). Furthermore, suppose Assumption \ref{ass:1} holds with $p\geq 6$ and $\hbeta$ converges at the rate in Theorem \ref{Supthm:1}. Then, for any $0<\delta <1$, there exists a constant $C_{(\eta_1, \mu_1, \|\tbeta\|_1, \|\tgamma\|_1)}$ such that the estimator \eqref{est:3} with $\tau = (\frac{nM_6}{\log(2d/\delta)})^{\frac{1}{6}}$, $\kappa = 14\sqrt{\frac{M_6\log(2d/\delta)}{n}}$, and $\lambda \geq C_{(\eta_1, \mu_1, \|\tbeta\|_1,\|\tgamma\|_1)}\kappa$ satisfies
	\begin{align*}
	P\bigg(\|\hgamma_5 - \ttgamma\|_2\leq 3\sqrt{s}\lambda\text{\ \ and\ \ } \|\hgamma_5 - \ttgamma\|_1\leq 12s\lambda\bigg)\geq 1-\delta.
	\end{align*}
	
\end{theorem}

Unlike the results in the previous sections, the argument in Theorem \ref{Sthm:3} holds for finite sample size $n$. The choice of the tuning parameter $\lambda$ depends on some quantities of $\tbeta$ and $\tgamma$, which need to be tuned in practice. Even though our estimation is built on the identity \eqref{equ:2}, the above result still holds for $\eta_1=0$.

The estimator based on the second-order Stein's identity is obtained similarly. Suppose the functions $f^2(\LD \xb, \tbeta\RD)$ and $g^2(\LD \xb, \tgamma\RD)$ together with $\xb$ satisfy the second-order regularity condition. Then
\begin{align}\label{equ:5}
\mE[y^2H(\xb)] = \mE [f^2(\LD \xb, \tbeta\RD)H(\xb)] + \mE [g^2(\LD \xb, \tgamma\RD)H(\xb)] = 2\eta_2\tbeta\tbetat + 2\mu_2\tgamma\tgammat,
\end{align}
where $\eta_2 = \mE[(f'(\xb^T\tbeta))^2]+ \mE[f(\xb^T\tbeta)f''(\xb^T\tbeta)]$ and $\mu_2$ is defined same as \eqref{equ:3}. Let $\tilde{U} =\frac{1}{n}\sum_{i=1}^{n}\Psi_{\tau}(y_i)^2\cdot\Psi_{\tau}(H(\xb_i))$ be the truncated counterpart of the left hand side in \eqref{equ:5}, $\barU = \ttU - \rbr{\hbeta^T\ttU\hbeta}\cdot\hbeta\hbeta^T$, and $\barV = \lamT(\barU, 1)$. Then our estimator can be computed as
\begin{align}\label{est:6}
\hgamma_6 \in \arg\max_{\|\bgamma\|_2\leq 1} |\bgamma^T\barV\bgamma|.
\end{align}

\begin{theorem}[Second-order orthogonal estimation]\label{Sthm:6}
	
	Suppose $f^2(\LD \xb, \tbeta\RD)$ and $g^2(\LD\xb, \tgamma\RD)$ together with $\xb$ satisfy the second-order regularity condition, $\|\tbeta\|_0 \vee \|\tgamma\|_0\leq s$, $\LD \tbeta,\tgamma\RD=0$, and $\mu_2\neq 0$ (defined in \eqref{equ:3}). Furthermore, suppose Assumption \ref{ass:1} holds with $p\geq 6$ and $\hbeta$ converges at the rate in Theorem \ref{Supthm:1}.	Then, for any $0<\delta<1$, there exists a constant $C'_{(\eta_2, \mu_2, \|\tbeta\|_1, \|\tgamma\|_1)}$ such that the estimator \eqref{est:6} with $\tau = (\frac{nM_6}{\log(2d^2/\delta)})^{\frac{1}{6}}$ and $\lambda\geq C'_{(\eta_2, \mu_2, \|\tbeta\|_1, \|\tgamma\|_1)}\sqrt{\frac{M_6\log(2d^2/\delta)}{n}}$ satisfies
	\begin{align*}
	P\bigg(\min_{\iota = \pm 1}\|\iota \hgamma_6 - \tgamma\|_2\leq \frac{4\sqrt{2}s\lambda}{\mu_2}\bigg)\geq 1-\delta.
	\end{align*}
	
\end{theorem}

We conclude this section by noting that while the requirement
$\tbetat\tgamma = 0$ might seem restrictive, our proof technique can be trivially modified to allow for a relaxed assumption stating that $\abr{\tbetat\tgamma} \lesssim \sqrt{\log d/n}$ for the same estimator, which would be often satisfied in a high-dimensional setting. Compared to the results in the last two subsections, here we utilize the approximate orthogonality to obtain non-asymptotic results, rather than relying on estimation of $f$. Furthermore, we note that the orthogonality condition holds in applications of generalized linear mixed models, where predictors that contribute to
the mean part will not be included in the variance part. Therefore, Theorem \ref{Sthm:3} and \ref{Sthm:6} are useful for orthogonal design generalized linear mixed models \citep{Faraway2016Extending, McCullagh2018Generalized}.

\section{Multiple Index Volatility Model}\label{sec:4}

In this section, we study the model \eqref{mod:2} with $v > 1$, which is called multiple index volatility model. We develop an estimator for $\tG$ based on the second-order Stein's identity. The first-order Stein's identity is not directly applicable here, unless combined with sliced inverse regression. See \cite{Babichev2018Slice, Lin2019Sparse} for related issues in multiple index models.

Our starting point is the second-order identity, which states that
\begin{align}\label{equ:6}
\mE[(y - f(\xb^T\tbeta))^2H(\xb)] = \mE[g^2(\tGt\xb)H(\xb)] = 2\tG \bLambda\tGt,
\end{align}
where $\bLambda = \mE[\nabla g(\tGt\xb)\nabla^Tg(\tGt\xb)+g(\tGt\xb)\nabla^2g(\tGt\xb)]\in\mR^{v\times v}$. Let $\mu_3 = \lambda_{\min}(\bLambda)$ be the minimum eigenvalue of $\bLambda$ and suppose $\mu_3>0$. Note that we could replace this identifiability condition by $\mu_3 = \lambda_{\max}(\bLambda)<0$. Our estimation procedure is similar to what we discussed in Section~\ref{sec:sivm:sec-order-estimation}, however, we will extract top $v$ eigenvectors that will estimate $\tG$ up to an orthogonal transformation.

In a low-dimensional setting, starting from $\hU$ defined in
\eqref{est:4}, we define $\hG_1$ as a solution to the following
optimization program:
\begin{equation}\label{Mest:1}
\begin{aligned}
\hG_1 \in \arg&\max_{\bG\in\mR^{d\times v}} \LD \hU, \bG\bG^T\RD\\
&\text{s.t. } \bG^T\bG = I_v.
\end{aligned}
\end{equation}

\begin{theorem}[Low-dimensional second-order estimator]\label{Mthm:1}
	
	Suppose conditions of Theorem \ref{Sthm:4} are satisfied and $\mu_3>0$. The estimator \eqref{Mest:1}, with $\hU$ set by Theorem \ref{Sthm:4}, satisfies
	\begin{align*}
	P\bigg(\inf_{Q\in\mQ^{v\times v}}\|\hG_1 - \tG Q\|_F \leq \frac{\Upsilon}{\mu_3}\big(d\sqrt{\frac{M_6\log(d/\delta)}{n}} + \frac{d\sqrt{M_6}\cdot \tte_{f, \delta}(n,1)}{\delta}\big)\bigg)\geq 1-\delta.
	\end{align*}
	Furthermore, if the mean estimation in Theorem \ref{Supthm:1} and \ref{Supthm:2} are satisfied, we have
	\begin{align*}
	\inf_{Q\in\mQ^{v\times v}}\|\hG_1 - \tG Q\|_F = O_P\rbr{ \frac{d}{\mu_3}\sqrt{\frac{\log d}{n}} }.
	\end{align*}
	
\end{theorem}

In a high-dimensional setting, we let $\hV = \mT_{\lambda}(\hU, v)$ be the top-$v$ sparse eigenvectors of $\hU$ where $\mT_{\lambda}(\cdot, \cdot)$ is defined in \eqref{equ:4}, then our high-dimensional estimator $\hG_2$ can be solved from \eqref{Mest:1} with $\hU$ replaced by $\hV$. Its statistical rate of convergence is given in next theorem.

\begin{theorem}[High-dimensional second-order estimator]\label{Mthm:2}
	
	Suppose conditions of Theorem \ref{Sthm:5} are satisfied (we replace $\|\tgamma\|_0\leq s$ by $\|\tG\|_{0,\max}\leq s$) and $\mu_3>0$. Under the same setup of $\lambda$ as in Theorem \ref{Sthm:5}, the estimator $\hG_2$ satisfies
	\begin{align*}
	P\bigg(\inf_{Q\in\mQ^{v\times v}}\|\hG_2 - \tG Q\|_F\leq \frac{4s\sqrt{v}\lambda}{\mu_3}\bigg)\geq 1-\delta.
	\end{align*}
	Furthermore, if the mean estimation in Theorem \ref{Supthm:1} and \ref{Supthm:2} are satisfied, we have
	\begin{align*}
	\inf_{Q\in\mQ^{v\times v}}\|\hG_2 - \tG Q\|_F = O_P\rbr{ \frac{s}{\mu_3}\sqrt{\frac{v\log d}{n}}  }.
	\end{align*}
\end{theorem}

Analogously, in the orthogonal case, i.e. $\tbetat\tG = \0$, we redefine $\barV$ in \eqref{est:6} by letting $\barV = \mT_{\lambda}(\barU, v)$, and apply \eqref{Mest:1} on $\barV$ to extract its top-$v$ eigenvectors. The estimator is denoted as $\hG_3$, and its convergence rate is summarized next.

\begin{theorem}[Orthogonal estimation]\label{Mthm:3}
	
	Suppose conditions of Theorem \ref{Sthm:6} are satisfied with $\tbetat\tG=\0$, $\|\tG\|_{0, \max}\leq s$, and $\mu_3>0$. Then, under the same setup of $\tau$ and $\lambda$ as in Theorem \ref{Sthm:6}, the estimator $\hG_3$ satisfies
	\begin{align*}
	\inf_{Q\in\mQ^{v\times v}}\|\hG_3 - \tG Q\|_F\leq \frac{4s\sqrt{v}\lambda}{\mu_3}
	\end{align*}
	with probability at least $1-\delta$.
	
\end{theorem}

The rate of convergence in the last two theorems is proportional to $s\sqrt{v}\lambda / \mu_3$ as $\tG\tGt$ has at most $s^2v$ nonzero elements when $\|\tG\|_{0,\max}\leq s$. If, instead, the sparsity structure on $\tG$ is assumed that $\|\tGt\|_{0, 0}\leq s$, i.e. $\tG$ has at most $s$ nonzero rows \citep{Xu2010Robust, Obozinski2011Supporta}, then $\tG\tGt$ has at most $s^2$ nonzero elements and the rate would be proportional to $s\lambda / \mu_3$.

In summary, the estimation of multiple index volatility models is a straight-forward generalization from the single index case. \cite{Dudeja2018Learning} considered a different problem: what if a multiple index model is misspecified to a single index model. We leave the problem of misspecification in the context of Stein's estimators for future work.

\section{Numerical Experiment}\label{sec:5}

We conduct extensive numerical experiments to validate theoretical results presented in Section \ref{sec:3} and \ref{sec:4}. We focus our attention on recovery of $\tG$ and verify convergence rates in a high-dimensional setting assuming the knowledge of $f$ and using  \cite{Yang2017High} to estimate $\tbeta$. Under this setup, all statements hold for finite sample. Also, only the step of estimating $\tbeta$ contributes the error $\tte_{f, \delta}(n, 1)$, which is still in smaller order comparing to the error occurs in the step of estimating $\tG$. Specifically, we will empirically show that the estimation error is upper bounded by $\sqrt{s\log d/n}$ for the first-order estimation and $s\sqrt{\log d/n}$ for the second-order estimation. The estimation accuracy is measured using (\ref{equ:cos}) for single index volatility models, and the sine distance, defined as
\begin{align*}
\|\sin(\angle(\hG, \tG))\|_F = \frac{1}{\sqrt{2}}\|\hG\hG^T - \tG\tGt\|_F,
\end{align*}
for multiple index volatility model. Throughout the simulations, we set the mean link function to be $f(x) = 2x+\cos(x)$ and consider three different designs for $\xb$: Gaussian, Student's $t$, and Gamma. Table \ref{tab:1} summarizes the parameters used in experiments, as well as the first- and second-order score functions. We let $\epsilon\sim \mN(0, 1)$. All simulation results are reported over 10 independent runs.

\begin{table}[h]
	\centering
	\begin{tabular}{ l | c |c | c }
		\hline
		Distribution & Parameter & First-order score & Second-order score \\
		\hline
		Gaussian & $\mu=0$, $\sigma = 1$ & $S(x) = x$ & $H(x) = x^2-1$ \\
		Student's $t$ & degree of freedom $13$ & $S(x) = \frac{14x}{13+x^2}$ & $H(x) = \frac{224x^2}{(13+x^2)^2} - \frac{14}{13+x^2}$\\
		Gamma & $k = 13$, $\theta = 2$ & $S(x) = \frac{1}{2} - \frac{12}{x}$ & $H(x) = \frac{132}{x^2} - \frac{12}{x} + \frac{1}{4}$ \\
		\hline
	\end{tabular}
	\caption{Distribution of covariate $\xb$}
	\label{tab:1}
\end{table}

Although our theorems rely on prior knowledge of score
	functions, we also consider the case in simulation that this
	information is not accessible. We apply the approach in
	\cite{Babichev2018Slice} to estimate score functions, which is
	originally inspired by the score matching method in
	\cite{Hyvaerinen2005Estimation}. In particular, suppose
	$\{\zeta_j(x)\}_{j=1}^m$ is a set of basis functions and the
	first-order score function satisfies
	$S(x) = \sum_{j=1}^{m}\nu_j^\star\zeta_j(x)$. We further define
	$\bzeta(x) = (\zeta_1(x); \ldots; \zeta_m(x))$ and
	$\bzeta'(x) = (\zeta_1'(x); \ldots; \zeta_m'(x))$. Then, given $n$
	independent samples, we estimate
	$\bnu^\star = (\nu_1^\star; \ldots; \nu_m^\star)$ by
	\begin{align}\label{equ:est:score}
	\hat{\bnu} = \rbr{ \frac{1}{n} \sum_{i=1}^{n} \bzeta(x_i)\bzeta(x_i)^T + \alpha I_m}^{-1}\rbr{\frac{1}{n}\sum_{i=1}^{n}\bzeta'(x_i)}.
	\end{align}
	Therefore, the first-order score function is estimated by
	$\hat{S}(x) = \hat{\bnu}^T\bzeta(x)$, and the second-order score function is estimated by
	$\hat{H}(x) = \hat{\bnu}^T\bzeta(x)\bzeta(x)^T\hat{\bnu} -
	\hat{\bnu}^T\bzeta'(x)$. In our simulations, we let $\alpha = 0.01$, $m = 101$ and, $\zeta_j(x) = \frac{1}{\sqrt{2\pi h^2}}\exp\rbr{-\frac{(x - z_j)^2}{2h^2}}$ is the Gaussian kernels with $h = 0.5$ and $z_j = -5 + 0.25(j-1)$, $j = 1, \ldots, m$.

\subsection{Single index experiment}\label{sec:8.1}

\begin{figure}[!t]
	\centering     
	\includegraphics[width=136mm]{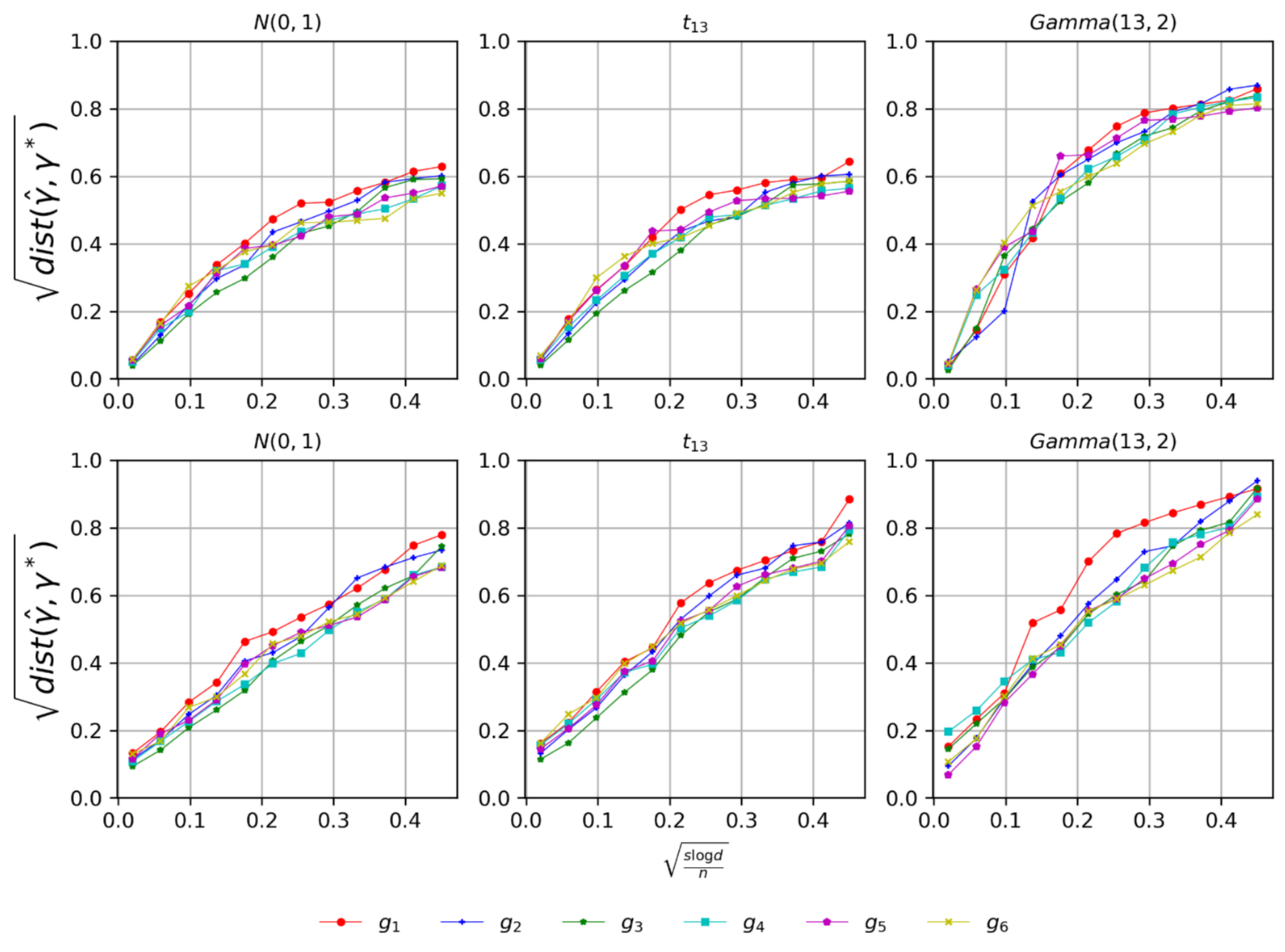}
	\caption{\textit{Estimation error for the first-order estimator $\hgamma_2$. The lines indicate different variance functions, while columns indicate different designs on $\xb$. The first row corresponds to the setup where score functions are known, while the bottom row to the setup where score functions are estimated by \eqref{equ:est:score}.}}\label{fig:1}
\end{figure}

\begin{figure}[!t]
	\centering     
	\includegraphics[width=136mm]{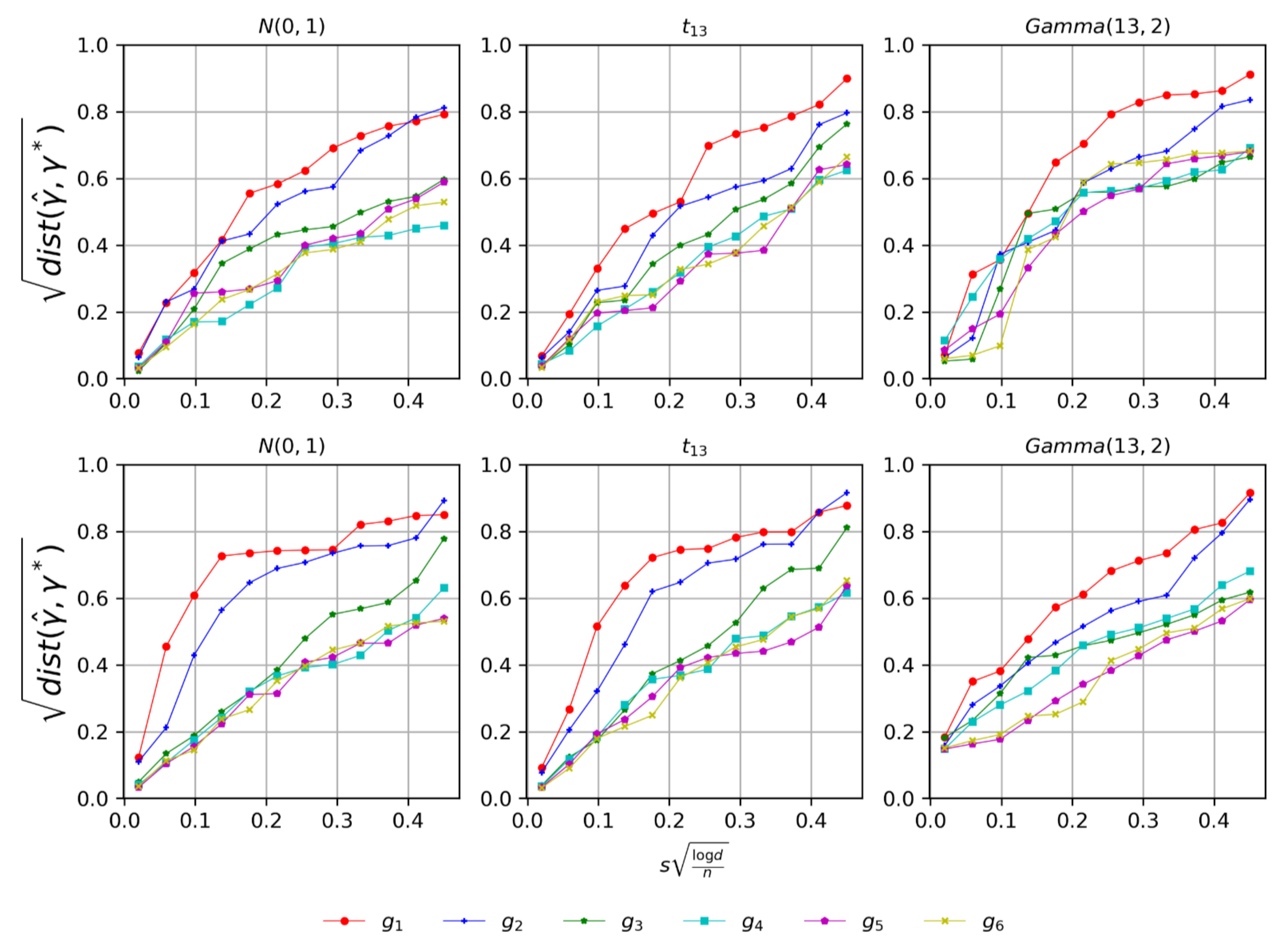}
	\caption{\textit{Estimation error for the second-order estimator $\hgamma_4$. Notations are same as Figure \ref{fig:1}.}}\label{fig:2}
\end{figure}

We consider two estimators of $\tgamma$ here: $\hgamma_2$ in \eqref{est:2} and $\hgamma_4$ in \eqref{est:5}. The optimization program in \eqref{equ:4} is approximately solved using the ADMM-based algorithm proposed in \cite{Vu2013Fantope} with  the Lagrange multiplier $\rho = 1$ (see equation (9) in \cite{Vu2013Fantope}). We consider six different variance link functions:
\begin{align*}
g_1(x) = x+1 + \cos(x); \quad g_2(x) = x + 1 + \exp(-x^2); \quad g_3(x) = x + 1 + \frac{\exp(x)}{(1 + \exp(x))^2};\\
g_4(x) = x^2 + x + \cos(x); \text{\ \ } g_5(x) = x^2+x+\exp(-x^2); \text{\ \ } g_6(x) = x^2+x + \frac{\exp(x)}{(1+\exp(x))^2}.
\end{align*}
We set $d = 100$ and $s = 10$, and vary the sample size $n$. For each $i\in[n]$ and $j\in[d]$, $[\xb_{i}]_j$ is generated independently from the corresponding distributions. For unknown coefficients $\tbeta$ and $\tgamma$, we first randomly generate positions of non-zero indices from $[d]$, and then each non-zero entry is set to $\pm 1/\sqrt{s}$ with equal probability. We set $\tau = 15(n/\log d)^{\frac{1}{6}}$ and $\lambda = 0.1\sqrt{\log d/n}$ for both first- and second-order estimations. The simulation results are shown in Figures \ref{fig:1} and \ref{fig:2}.

From the plots, we observe that our theories explain the relationship between the estimation error and the problem parameters. In particular, when sample size is sufficiently large, Figure \ref{fig:1} shows that the estimation error for the first-order method decreases linearly with $\sqrt{s\log d/n}$, as suggested by Theorem~\ref{Sthm:2}. Although there are small differences between different designs, we observe that for each design the error decreases linearly in the control parameters. Furthermore, even if score functions are unknown in practice, we note from the second row of Figure \ref{fig:1} that one can estimate score functions first and then plug estimated functions into our procedure. The resulted error is still controllable when sample size is large. Similar observations hold for the estimator $\hgamma_4$ which is based on the second-order Stein's identity. We should also mention that, for the second-order method, the performance of our  procedure on linear functions $\{g_i\}_{i=1}^3$ is worse than quadratic functions $\{g_i\}_{i=4}^6$, since their second derivative contains less information of signals.

\subsection{Multiple index experiment}

\begin{figure}[!t]
	\centering     
	\includegraphics[width=136mm]{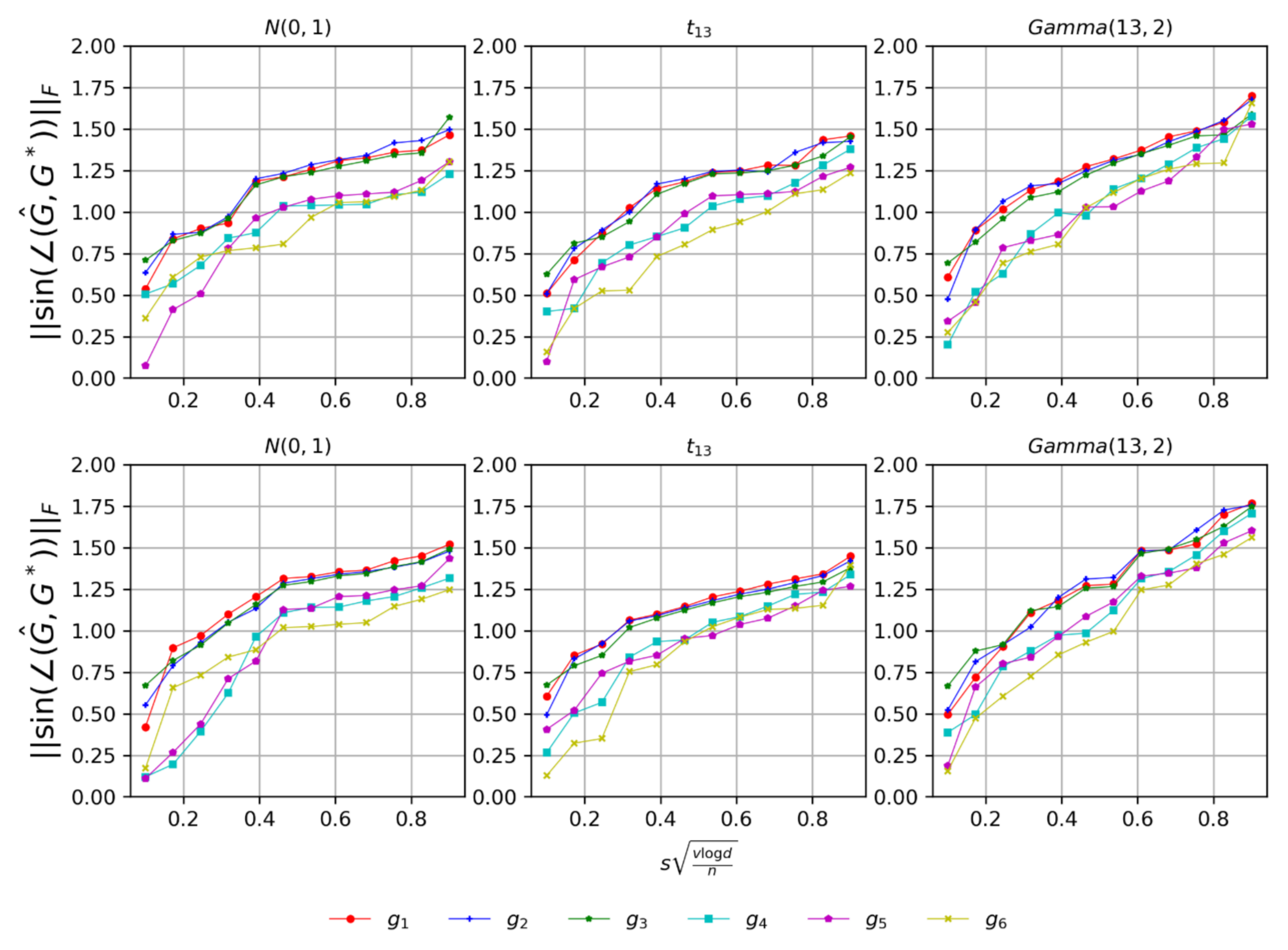}
	\caption{\textit{Estimation error for the second-order estimator $\hG_2$ (regular case). }}\label{fig:3}
\end{figure}

\begin{figure}[!t]
	\centering     
	\includegraphics[width=136mm]{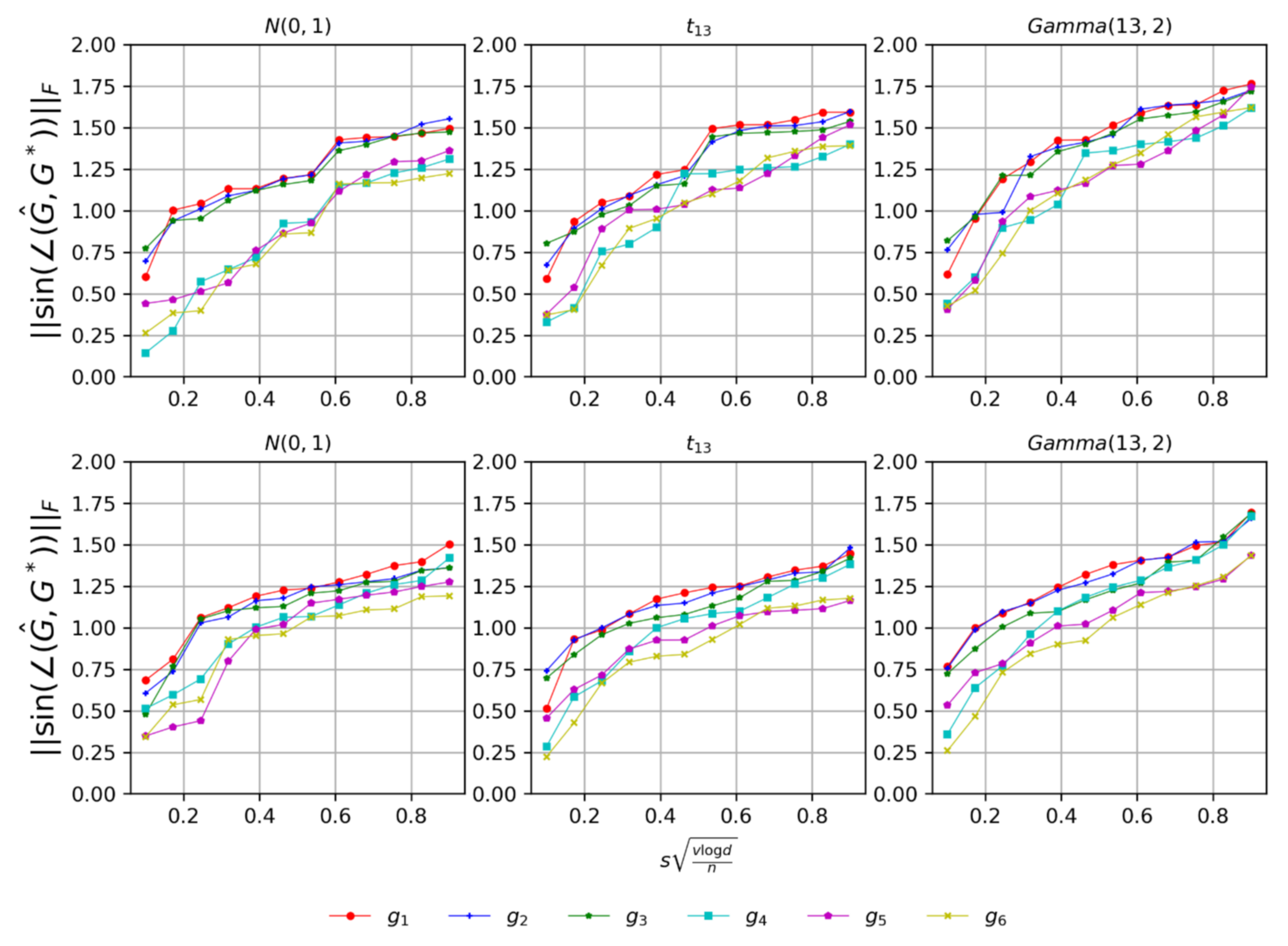}
	\caption{\textit{Estimation error for the second-order estimator $\hG_3$ (orthogonal case). }}\label{fig:4}
\end{figure}

We consider two estimators of $\tG$ here: $\hG_2$ and $\hG_3$ described in Section~\ref{sec:4}. We set $d = 100$, $s = 10$, and $v=3$. We let $\text{supp}(\tgamma_j) = [(j-1)s, js]$ and each entry in the support is set to $\pm 1/\sqrt{s}$ with equal probability. The variance link function is defined as $g(\tGt\xb) = \sum_{j=1}^vg(\xb^T\tgamma_j)$ where $g\in\{g_i\}_{i=1}^6$ is the variance function used in the single index experiment. In this experiment, we let $\tau = 10(n/\log d)^{\frac{1}{6}}$, $\lambda = 0.01\sqrt{\log d/n}$. The results are shown in Figures~\ref{fig:3} and \ref{fig:4}, and again we observe that the error rate is correctly explained by our theory for sufficiently large sample size. In particular, we observe that the error rate is sublinear in term of $s\sqrt{v\log d/n}$ as expected.

\section{Discussion}\label{sec:6}

We proposed new estimators for parametric components of index volatility models based on the first- and second-order Stein's identities. Our approach extracts the direction of variance index by multiplying the score variables with residuals and incorporating weighted mean squared minimization with truncation and regularization to accommodate for estimation in a high-dimensional and heavy-tailed setting. We established statistical rates of convergence of our estimators under both single and multiple index structures, which were then verified in finite samples through numerical experiments. In particular, our results were qualitatively the same under a range of designs, including heavy-tailed ones. The estimation rate is comprised
of both nonparametric and parametric rates, though the parametric rate is the dominant term and matches the corresponding rate in the mean estimation. We illustrated that when the mean index is orthogonal to the variance index, which would naturally be the case in many high-dimensional applications, we need not estimate the mean link function and can obtain finite sample results.

Our estimation procedures rely on estimation of the gradient of link function through Stein's identities, which is aligned with the
e.d.r.~direction in index volatility models. Using this approach, we were able to naturally extend traditional fixed design setup to
randomized design. The drawback of the approach based on Stein's identity is that the prior knowledge of a distribution of covariates $\xb$ is needed. In a low-dimensional setting,
\cite{Babichev2018Slice} proposed an approach for estimating the first-order score function under the assumption that the score function can be represented as a finite linear combination of basis functions in a given dictionary. We implement this approach in our experiments to address practical concerns, however theoretical analysis for extending their approach to a high-dimensional setting seems challenging without strong assumptions on the underlying distribution of $\xb$. We leave investigation of possible estimators for future work. Fortunately, in some applications, such as compressed sensing \citep{Ai2014One, Davenport20141} or phase retrieval
\citep{Candes2015Phasea, Candes2015Phase}, the distribution of
covariates is known.

In this work, we have focused on point estimation of the e.d.r.~direction. Establishing tools that would allow for construction of confidence intervals and more generally uncertainty quantification is an important future direction. Another direction is improving the statistical rate obtained by the second-order estimation procedures, which arises from the finite moment condition on $H(\xb)$ that makes bounding of the restricted operator norm difficult. We note, however, that obtaining the error rate of $\sqrt{s\log d/n}$ under heavy-tailed
design for sparse PCA is still an open problem. In addition, developing a robust version of the estimators based on absolute
residuals, as used for robust variance estimation \citep{Davidian1987Variance}, would further enlarge potential for diverse applications of our methodology. Finally, it is also of interest to study higher-order methods. \cite{Dudeja2018Learning} proposed $l$-th order estimator based on Hermite polynomials, and proved that the rate of convergence is $\sqrt{d(\log d)^l/n}$. This is useful especially when $\mu_1 = \mu_2= 0$. Analogously and to deal with non-Gaussian design, one can recursively define higher-order score variable $S_l(\xb) = - S_{l-1}(\xb)\otimes \nabla \log p(\xb) + \nabla S_{l-1}(\xb)$ and use higher-order Stein's identity, $\mE[f(\xb^T\tbeta) S_l(\xb)] = f^{(l)}(\xb^T\tbeta)\cdot (\tbeta)^{\otimes l}$, to estimate the index. However, this may require tensor decomposition and is less efficient than the nonconvex method in \cite{Dudeja2018Learning}. How to bridge this computational gap and have an efficient higher-order estimator with optimal convergence rate is an attractive topic for future work.

\section*{Acknowledgments}

This work is partially supported by the William S. Fishman Faculty Research Fund at the University of Chicago Booth School of Business. This work was completed in part with resources provided by the University of Chicago Research Computing Center.

\bibliographystyle{my-plainnat}
\bibliography{paper}

\appendix
\pagebreak
\numberwithin{equation}{section}

\begin{center}
	\textbf{\large Supplemental Materials: \\ High-dimensional Index Volatility Models via Stein's Identity}
\end{center}

\section{Estimation of Index Mean Function}\label{sec:A}

In this section, we present results on the mean estimation for the model \eqref{mod:2}. In particular, we develop an explicit formula for $\bare_f(\hbeta^T\xb,n,1)$ in \eqref{cond:5}, and further derive a bound on its first moment, $\hate_f(\hbeta, n, 1)$, and error rate $\tte_{f, \delta}(n, 1)$. By Corollary \ref{cor:1} we finally verify the argument in \eqref{cond:7}. Our estimation procedure is based on two steps. First, we use the approach in \cite{Yang2017High} to estimate the mean index $\tbeta$. Next, a local linear regression is applied on the pair $(y, \hbeta^T\xb)$ to obtain the mean link function estimator $\hf$. Finally, we use $\hf(\LD \hbeta,\xb\RD)$ as an estimation of $f(\LD \tbeta, \xb\RD)$. To simplify the analysis, we assume that two steps are conducted on independent samples of size $n$ for each, which is obtained, for example, by sample splitting in practice. This simplifies the analysis while keeping the statistical rate unchanged. We note that the local linear regression is just one way to estimate the nonparametric component in index models. See \cite{Liu2013robust} for a robust estimator as an alternative.

To unify the presentation, we consider a more general heteroscedastic index model:
\begin{align}\label{mod:3}
y\mid \xb = f(\LD\tbeta, \xb\RD) + \ttg(\xb)\epsilon,
\end{align}
where $\mE[\epsilon | \xb] = 0$, $\mE[\epsilon^2 | \xb] = 1$. By setting $\ttg(\xb) = g(\LD \tgamma,\xb\RD)$ we obtain the single index volatility model and $\ttg(\xb) = g(\tGt\xb)$ would lead to multiple index volatility model. The estimator $\hbeta$ in \cite{Yang2017High} is defined as
\begin{align}\label{est:beta}
\hbeta = \phi_{\lambda}\bigg(\frac{1}{n}\sum_{i=1}^n\Psi_{\tau}(y_i)\Psi_{\tau}(S(\xb_i))\bigg),
\end{align}
where $\lambda$ and $\tau$ are tuning parameters. This Lasso-type estimator comes from the first-order Stein's identity applied on the response $y$. Here $\tau$ is the truncation threshold and $\lambda$ controls the sparsity of $\hbeta$. Note that $\hbeta$ can be computed without the knowledge of $f$. Its convergence rate is presented in the following theorem.

\begin{theorem}[$\tbeta$ estimation]\label{Supthm:1}
	
Suppose the function $f(\LD\tbeta,\xb\RD)$ together with $\xb$ satisfies the first-order regularity condition. Furthermore, suppose Assumption \ref{ass:1} holds with $p \geq 4$ and $\|\tbeta\|_0\leq s$. Then, for any $0<\delta<1$, the estimator $\hbeta$ defined in \eqref{est:beta} for $\tbeta$ in model \eqref{mod:3}, with $\lambda \asymp \sqrt{\log(d/\delta)/n}$ and $\tau\asymp (n/\log(d/\delta))^{1/4}$, satisfies
\begin{align*}
\|\hbeta - \tbeta\|_2\lesssim\sqrt{\frac{s\log(d/\delta)}{n}},\quad \|\hbeta - \tbeta\|_1\lesssim s\sqrt{\frac{\log(d/\delta)}{n}}, \quad \|\hbeta - \tbeta\|_{\infty}\lesssim\sqrt{\frac{\log(d/\delta)}{n}},
\end{align*}
with probability at least $1-\delta$.

\end{theorem}

Theorem 4.2 in \cite{Yang2017High} proves the above result for a high-dimensional homoscedastic single index model. Theorem~\ref{Supthm:1} states a more general result that is valid for a heteroscedastic single index model. The proof is omitted as the proof strategy of \cite{Yang2017High} is directly applicable, since 
\begin{align*}
\mE[yS(\xb)] = \mE[f(\LD \tbeta, \xb\RD)S(\xb)] + \mE[\epsilon\ttg(\xb)S(\xb)] = \mE[f(\LD \tbeta, \xb\RD)S(\xb)] = c\tbeta.
\end{align*}

With $\hbeta$ defined above, we use a local linear regression estimator to estimate $f(\cdot)$. The results borrows from \cite{Fan1993Local} and requires the following assumption (see Condition 1 in \cite{Fan1993Local} for comparison). The assumption is specifically required for estimation using local linear regression and would need to be adjusted if other estimators are used.

\begin{assumption}\label{ass:3}
	
For any fixed estimator $\hbeta$ with the rate of convergence as in Theorem \ref{Supthm:1}, we assume
\begin{enumerate}[label=(\alph*),topsep=1pt]
\setlength\itemsep{0.3em}
\item (smoothness) $f\in C^2(\mR)$ with $|f'(x)|\vee |f''(x)|\leq L_1$, $\forall x\in\mR$.
\item ($4$th moment function) $\rbeta(x) =\mE[\ttg^4(\xb) \mid \hbeta^T\xb = x]$ is continuous and bounded.
\item (density) $\hbeta^T\xb$ has bounded positive density $q_{\hbeta}$ with $|q_{\hbeta}(x) - \qbeta(y)|\leq L_2|x-y|^{\alpha}$ for $\alpha\in(0,1)$.
\item (tail) there exists a constant $L_3>0$ such that $|f(x)|^4/\qbeta^8(x)$, $\rbeta(x)/\qbeta(x)$, and $|f'(x)|^4/\qbeta(x)$ are integrable on $(-\infty, -L_3)\cup(L_3,\infty)$\footnote{Here we assume $\qbeta$ has support over $\mR$. If not, functions only need to be integrable on the tail of the support.}.
\end{enumerate}
	
\end{assumption}

With these assumptions, we define the local linear estimator as
\begin{align}\label{est:f}
\hf(t) = \frac{\sum_{i=1}^n w_iy_i}{\sum_{i=1}^nw_i + (nh)^{-2}},
\end{align}
where $w_i = K_h(\xb_i^T\hbeta - t)\big(s_{n2} - (\xb_i^T\hbeta - t)\cdot s_{n1}\big)$ with $s_{nl} = \sum_{i=1}^n K_h(\xb_i^T\hbeta - t)(\xb_i^T\hbeta - t)^l$ for $l = 1,2$, $h$ is the bandwidth, and 
$K_h(\cdot) = K(\cdot/h)/h$ for a kernel function $K(\cdot)$. With a Gaussian kernel, we have the following rate of convergence for the mean estimator.

\begin{theorem}[Mean estimation]\label{Supthm:2}
	
Suppose $\hbeta$ is estimated using samples $\D_1$ and satisfies the rate in Theorem \ref{Supthm:1}. Given $\hbeta$, let $\hf$, defined in \eqref{est:f}, be a local linear estimator of the link function $f$	using samples $\D_2 = \{y_i, \hbeta^T\xb_i\}_{i\in[n]}$ that are independent from $\D_1$. Suppose Assumption~\ref{ass:1} with $p\geq 4$ and Assumption \ref{ass:3} (a-c) hold and the bandwidth $h$ satisfies $h\rightarrow 0$ and $nh\rightarrow \infty$. Then there exist $N$ and $\Upsilon_1$ such that $\forall n\geq N$,
\begin{align}\label{a:2}
\mE[|\hf(\xb^T\hbeta) - f(\xb^T\tbeta)|^4 \mid \xb, \D_1]\leq  \Upsilon_1\bigg(h^8 +& \frac{\rbeta(\xb^T\hbeta)}{n^2h^2\qbeta^2(\xb^T\hbeta)} + \frac{\|\hbeta - \tbeta\|_2^4\cdot|f'(\xb^T\hbeta)|^4}{h^2\qbeta^2(\xb^T\hbeta)} \nonumber \\
+& |\xb^T\tbeta - \xb^T\hbeta|^4 + \frac{|f(\xb^T\hbeta)|^4}{n^{16}h^{16}\qbeta^8(\xb^T\hbeta)}\bigg).
\end{align}
Furthermore, suppose Assumption \ref{ass:3} (d) holds as well and $h\asymp n^{-1/5}$, then there exists a constant $\Upsilon_2$ such that
\begin{equation}\label{a:3}
P\bigg(\biggl\|\mE\big[|\hf(\xb^T\hbeta) - f(\xb^T\tbeta)|^2\cdot W(\xb) \mid \D_1, \D_2\big]\biggr\|_{\infty}\geq \frac{\sqrt{M_2}}{\delta}\underbrace{\Upsilon_2n^{-4/5}}_{\tte_{f, \delta}(n, 1)}\bigg)\leq 2\delta,
\end{equation}
where $W(\xb)$ is either the first-order score variable $S(\xb)$ or the second-order score variable $H(\xb)$. 
	
\end{theorem}

We see equation \eqref{a:2} provides an explicit form of $\bare_f(\hbeta^T\xb, n, 1)$ in \eqref{cond:5}. We explicitly write out higher-order terms in \eqref{a:2} to clarify the difference between a high-dimensional single index model and a nonparametric model. In a low-dimensional setting with $\xb$ being in a compact set and $\qbeta$ lower bounded away from zero \citep{Zhu2006Empirical, VanKeilegom2010Semiparametric, Wang2010Estimation, Lian2015Variance}, the last two terms can be ignored.

As discussed in Section \ref{sec:2}, estimation of $\tG$ is possible under a heavy-tailed design if condition \eqref{cond:6} holds. Assumption \ref{ass:3} (d) implies that expectation of the right hand side of \eqref{a:2}, conditional on $\hbeta$, is bounded on the tails; while within the interval $[-L_3, L_3]$, we can make use of continuity so that the integral is bounded naturally. In particular, the assumption imposes conditions on the decay rate of $|f(x)|$, $|f'(x)|$, and $\rbeta(x)$, and holds for any random variables that have compact support. Taking conditional expectation for $\bare_f(\hbeta^T\xb, n, 1)$ and ignoring all smaller order terms, we have 
\begin{align*}
\hate_{f}(\hbeta, n, 1)\lesssim h^8 + \frac{1}{n^2h^2} + \frac{\|\hbeta - \tbeta\|_2^4}{h^2}.
\end{align*}
Moreover, using the fact that $\|\hbeta - \tbeta\|_2\lesssim 1/\sqrt{n}$, we can set the bandwidth $h\asymp n^{-1/5}$ to obtain $\tte_{f, \delta}(n, 1)\lesssim n^{-4/5}$.

We also point out that $L_2$ and $L_3$ in Assumption \ref{ass:3} (c-d) do not depend on $\hbeta$ as long as it is close to $\tbeta$, which is assumed similarly in \cite{Zhang2018Quasi}. An equivalent statement would be to assume $q_{\tbeta}, r_{\tbeta}$ satisfy conditions (c-d) and further add some continuity conditions on $q_{\tbeta}, r_{\tbeta}$ with respect to $\bbeta$, such that $|q_{\hbeta} - q_{\tbeta}|$ and $|r_{\hbeta} - r_{\tbeta}|$ are small enough. For example, suppose $|\qbeta(x) - \qbeta(x)|\leq L_x\|\hbeta - \tbeta\|_2$ and $\sup_x L_x<\infty$, then by triangle inequality we have $\forall x, y$
\begin{align*}
|\qbeta(x) - \qbeta(y)|\leq& |\qbeta(x) - q_{\tbeta}(x)| + |q_{\tbeta}(x) + q_{\tbeta}(y)| + |\qbeta(y) - q_{\tbeta}(y)|\\
\leq & \underbrace{2\sup_x L_x\|\hbeta - \tbeta\|_2 }_{\text{ignorable}}+ L_3|x- y|^{\alpha} \lesssim L_3'|x - y|^{\alpha}.
\end{align*}

\section{Proofs of Main Theorems and Lemmas}\label{sec:B}

Throughout this section, we write $M$, omitting the subscript of $M_p$, since the moment degree $p$ is clear from the statement of theorem. For simplicity, we replace the truncation function $\Psi_{\tau}$ by notation $\widecheck{(\cdot)}$ where the truncation threshold $\tau$ is contained implicitly. In particular, we have $\Psi_{\tau}(\bv) = \widecheck{\bv}$. We use $\Upsilon>0$ to denote a generic constant, which may take different values for each appearance. For Theorem \ref{Sthm:1}, \ref{Sthm:2}, \ref{Sthm:4}, \ref{Sthm:5}, \ref{Mthm:1}, \ref{Mthm:2}, we only prove the first part of statement since the second part is trivial to obtain by plugging in $\tte_{f, \delta}(n, 1)\asymp n^{-4/5}$.

\subsection{Proof of Theorem \ref{thm:mean}}

We prove the result for $W(\xb) = S(\xb)$ as the other case is shown analogously. For any $j\in[d]$, 
\begin{align*}
\mEf\big[|\hf(\xb) - f(\xb)|^2\cdot|S(\xb)_j|\big]\leq& \sqrt{\mEf\big[|\hf(\xb) - f(\xb)|^4\big]}\sqrt{\mE[|S(\xb)_j|^2]} \\
\leq &\sqrt{M}\sqrt{\mEf\big[|\hf(\xb) - f(\xb)|^4\big]}.
\end{align*}
Therefore, we have
\begin{align*}
\big\|\mEf\big[|\hf(\xb) - f(\xb)|^2\cdot S(\xb)\big]\big\|_{\infty}\leq \sqrt{M}\sqrt{\mEf\big[|\hf(\xb) - f(\xb)|^4\big]}.
\end{align*}
By Markov's inequality, for any $0<\delta<1$, with probability $1-\delta$, we have
\begin{align*}
\sqrt{\mEf\big[|\hf(\xb) - f(\xb)|^4\big]}\leq \sqrt{\frac{\mE[|\hf(\xb) - f(\xb)|^4]}{\delta}}\stackrel{(\ref{cond:3})}{\leq} \frac{\tte_f(n, d)}{\delta}.
\end{align*}
Combining the above two inequalities completes the proof.

\subsection{Proof of Corollary \ref{cor:1}}

Similar to the proof of Theorem \ref{thm:mean}, we have for any $j\in[d]$, 
\begin{align*}
\mEbf\big[|\hf(\xb^T\hbeta) - f(\xb^T\tbeta)|^2\cdot|S(\xb)_j|\big]\leq & \sqrt{\mEbf\big[|\hf(\xb^T\hbeta) - f(\xb^T\tbeta)|^4\big]}\sqrt{\mE[|S(\xb)_j|^2]} \\
\leq &\sqrt{M}\sqrt{\mEbf\big[|\hf(\xb^T\hbeta) - f(\xb^T\tbeta)|^4\big]}
\end{align*}
and
\begin{align*}
\big\|\mEbf\big[|\hf(\xb^T\hbeta) - f(\xb^T\tbeta)|^2\cdot|S(\xb)_j|\big]\big\|_{\infty}\leq \sqrt{M}\sqrt{\mEbf\big[|\hf(\xb^T\hbeta) - f(\xb^T\tbeta)|^4\big]}.
\end{align*}
By Markov's inequality, for any $0<\delta<1$ and any sample set $\D_1$,
\begin{align}\label{C:1}
P\bigg(\mEbf\big[|\hf(\xb^T\hbeta) - f(\xb^T\tbeta)|^4\big]\leq \frac{\mEb\big[|\hf(\xb^T\hbeta) - f(\xb^T\tbeta)|^4\big]}{\delta} \Mid  \D_1\bigg)\geq 1 - \delta,
\end{align}
where the probability is taken over randomness in $\D_2$. By the definition in \eqref{cond:5}, we have
\begin{align*}
P\bigg(\sqrt{\mEbf\big[|\hf(\xb^T\hbeta) - f(\xb^T\tbeta)|^4\big]}\leq \sqrt{\hate_f(\hbeta, n, 1)}/\delta \Mid  \D_1\bigg)\geq 1 - \delta.
\end{align*}
Under the condition \eqref{cond:6}, we have
\begin{align*}
P\bigg(\bigg\|&\mEbf\big[|\hf(\xb^T\hbeta) - f(\xb^T\tbeta)|^2\cdot|S(\xb)_j|\big]\bigg\|_{\infty}\geq \frac{\sqrt{M}\cdot\tte_{f, \delta}(n, 1)}{\delta}\bigg)\\
\leq & P\bigg(\sqrt{\mEbf\big[|\hf(\xb^T\hbeta) - f(\xb^T\tbeta)|^4\big]} \geq \tte_{f, \delta}(n, 1)/\delta\bigg)\\
\leq & P\bigg(\sqrt{\mEbf\big[|\hf(\xb^T\hbeta) - f(\xb^T\tbeta)|^4\big]}\geq \sqrt{\hate_{f}(\hbeta, n, 1)}/\delta \text{\ or\ }  \sqrt{\hate_f(\hbeta, n, 1)}\geq \tte_{f, \delta}(n, 1)\bigg)\\
\leq &P\bigg(\sqrt{\mEbf\big[|\hf(\xb^T\hbeta) - f(\xb^T\tbeta)|^4\big]}\geq \sqrt{\hate_f(\hbeta, n, 1)}/\delta\bigg) + P\bigg(\sqrt{\hate_f(\hbeta, n, 1)}\geq \tte_{f, \delta}(n, 1)\bigg)\\
\stackrel{\eqref{cond:6}}{\leq} & \delta + \int P\bigg(\sqrt{\mEbf\big[|\hf(\xb^T\hbeta) - f(\xb^T\tbeta)|^4\big]}\geq \sqrt{\hate_f(\hbeta, n, 1)}/\delta\Mid \D_1\bigg)\underbrace{dP(\D_1)}_{\substack{\text{take integral over}\\ \\ \text{randomness in } \D_1}}\\
\leq& 2\delta.
\end{align*}
Here the last inequality is due to the fact that $\D_1$ and $\D_2$ are independent, so \eqref{C:1} holds uniformly for any $\D_1$.

\subsection{Proof of Theorem \ref{Sthm:1}}

Since the samples used for estimating $\ttgamma$ are independent from $\hf$, $\hbeta$,
\begin{align}\label{b:12}
\mEbf[(y - \hf(\xb^T\hbeta))^2\cdot S(\xb)] & = \mEbf [(y - f(\xb^T\tbeta) + f(\xb^T\tbeta) - \hf(\xb^T\hbeta))^2\cdot S(\xb)] \nonumber\\
& = \mE[g^2(\xb^T\tgamma)S(\xb)] + \mEbf[(f(\xb^T\tbeta) - \hf(\xb^T\hbeta))^2\cdot S(\xb)] \nonumber\\
& \stackrel{\eqref{equ:1}}{=} \ttgamma + \mEbf[(f(\xb^T\tbeta) - \hf(\xb^T\hbeta))^2\cdot S(\xb)].
\end{align}
For the second term, according to \eqref{cond:7}, for any $\delta>0$
\begin{align}\label{b:13}
P\bigg(\biggl\|\mEbf\big[|\hf(\LD\hbeta, \xb\RD) - f(\LD\tbeta, \xb\RD)|^2\cdot S(\xb)\big]\biggr\|_{\infty}\geq \frac{\sqrt{M}\cdot\tte_{f, \delta}(n,1)}{\delta} \bigg)\leq 2\delta.
\end{align}
(Here $M$ refers to $M_6$ since $\sqrt{M_2}\leq \sqrt{M_6}$.) Next, we bound the error that occurs when using $\hgamma_1$ to approximate the left hand side term in \eqref{b:12}. We apply Lemma \ref{aux1:lem:2}. Based on Assumption \ref{ass:1} ($p\geq 6$) we know that for some constant $\Upsilon_1$
\begin{align*}
\mE[|y|^6] \vee \mE[|S(\xb)_j|^6] \vee \mE[|f(\tbetat\xb)|^6]\leq \Upsilon_1 M, \text{\ \ } \forall j\in[d]. 
\end{align*}
Furthermore, by Markov's inequality, for any $\delta>0$, we have
\begin{align*}
\mEbf[|\hf(\xb^T\hbeta)|^6]\leq & \frac{32}{\delta} \mE[|\hf(\xb^T\hbeta) - f(\xb^T\hbeta)|^6]+ 32\mEb[|f(\xb^T\hbeta)|^6]
\end{align*}
with probability $1-\delta$. The first term goes to zero as $n\rightarrow \infty$ and it only attributes to higher order errors. Thus, there exists $N_{\delta}$ such that $\forall n\geq N_{\delta}$
\begin{align*}
\frac{32\mE[|\hf(\xb^T\hbeta) - f(\xb^T\hbeta)|^6]}{\delta} \leq \frac{\Upsilon_1M}{2}.
\end{align*}
Roughly, by nonparametric rate, we only require $\frac{1}{\delta n^p}\leq M$, which implies $n\geq (\frac{1}{\delta M})^{\frac{1}{p}}\eqqcolon N_{\delta}$ for some $p>1$. From Lemma \ref{aux1:lem:1}, we also have $\mEb[|f(\xb^T\hbeta)|^6]\leq \frac{\Upsilon_1M}{64}$ for a sufficiently large $N$ 
(not depending on $\delta$) and $\Upsilon_1$. Combining them together, we have $\mEbf[|\hf(\xb^T\hbeta)|^6]\leq \Upsilon_1M$ with probability $1-\delta$. Based on the definition of $\hgamma_1$ in \eqref{est:1},
\begin{align*}
\|\hgamma_1 -& \mEbf[(y - \hf(\xb^T\hbeta))^2\cdot S(\xb)]\|_{\infty}\\
\leq & \biggl\|\frac{1}{n}\sum_{i=1}^{n}\wyi^2\wSXi - \mEbf[y^2S(\xb)]\biggr\|_{\infty} + 2\biggl\|\frac{1}{n}\sum_{i=1}^{n}\wyi\wfxi\wSXi - \mEbf[y\hf(\xb^T\hbeta)S(\xb)]\biggr\|_{\infty}\\
& + \biggl\|\frac{1}{n}\sum_{i=1}^{n}\wfxi^2\wSXi - \mEbf[\hf(\xb^T\hbeta)S(\xb)]\biggr\|_{\infty}.
\end{align*}
For the above each term, applying Lemma \ref{aux1:lem:2}, setting $\tau= (\frac{n\Upsilon_1 M}{\log(2d/\delta)})^{\frac{1}{6}}$ and taking union bound over all $d$ indices, we obtain
\begin{align}\label{b:14}
P\bigg(\biggr\|\hgamma_1 - \mEbf[(y - \hf(\xb^T\hbeta))^2S(\xb)]\biggl\|_{\infty}\leq 28\sqrt{\frac{\Upsilon_1 M\log(2d/\delta)}{n}}\bigg)\geq 1-4\delta.
\end{align}
Combining \eqref{b:12}, \eqref{b:13} and \eqref{b:14} together and 
adjusting $\delta$ properly, we know there exists constant $\Upsilon$ large enough such that, by setting $\tau= \Upsilon(\frac{n M}{\log(d/\delta)})^{\frac{1}{6}}$,
\begin{align}\label{b:15}
\|\hgamma_1 - \ttgamma\|_{\infty}\leq & \|\hgamma_1 - \mEbf[(y - \hf(\xb^T\hbeta))^2S(\xb)]\|_{\infty} + \|\mEbf[(f(\xb^T\tbeta) - \hf(\xb^T\hbeta))^2S(\xb)]\|_{\infty} \nonumber\\
\leq & \Upsilon\bigg(\sqrt{\frac{M\log(d/\delta)}{n}} + \frac{\sqrt{M}\cdot\tte_{f, \delta}(n,1)}{\delta}\bigg)
\end{align}
with probability at least $1-\delta$. This completes the first part of the proof.

Moreover, since $\|\hgamma_1 - \ttgamma\|_{2}\leq \sqrt{d}\|\hgamma_1 - \ttgamma\|_{\infty}$, we know that
\begin{align}\label{b:16}
P\bigg(\|\hgamma_1 - \ttgamma\|_{2}\leq \Upsilon\big(\sqrt{\frac{Md\log(12d/\delta)}{n}} + \frac{\sqrt{Md}\cdot\tte_{f, \delta}(n,1)}{\delta}\big)\bigg)\geq 1-\delta.
\end{align}
Plugging in $\tte_{f, \delta}(n, 1)\asymp n^{-4/5}$, we see that the latter rate is negligible. Thus, it suffices to show $\dist(\hgamma_1 , \tgamma)\lesssim \frac{1}{\mu_1^2}\|\hgamma_1 - \ttgamma\|_2^2$. In fact, because $\|\tgamma\|_2 = 1$,
\begin{align}\label{b:17}
\dist(\hgamma_1, \tgamma) = &1 - \frac{|\hgamma_1^T\tgamma|}{\|\hgamma_1\|_2} = 1 - \frac{1}{2|\mu_1|\|\hgamma_1\|_2}\big|\|\ttgamma\|_2^2 + \|\hgamma_1\|_2^2 - \|\hgamma_1 - \ttgamma\|_2^2\big| \nonumber\\
\leq & \bigg(1 - \frac{|\mu_1|}{2\|\hgamma_1\|_2} - \frac{\|\hgamma_1\|_2}{2|\mu_1|}\bigg) + \frac{1}{2|\mu_1|\|\hgamma_1\|_2}\cdot\|\hgamma_1 - \ttgamma\|_2^2\nonumber\\
\leq &\frac{1}{2|\mu_1|^2}\cdot\|\hgamma_1 - \ttgamma\|_2^2\cdot\frac{|\mu_1|}{\|\hgamma_1\|_2}.
\end{align}
Note that
\begin{align*}
|\mu_1| - \|\hgamma_1 - \ttgamma\|_2\leq \|\hgamma_1\|_2\leq |\mu_1| + \|\hgamma_1 - \ttgamma\|_2,
\end{align*}
hence we have ($\|\hgamma_1 - \ttgamma\|_2/|\mu_1|\leq 1/2$)
\begin{align*}
1 - \frac{\|\hgamma_1 - \ttgamma\|_2}{|\mu_1|}\leq& \frac{|\mu_1|}{\|\hgamma_1\|}\leq 1+\frac{2\|\hgamma_1 - \ttgamma\|_2}{|\mu_1|}.
\end{align*}
Plugging back into \eqref{b:17} concludes the proof.

\subsection{Proof of Theorem \ref{Sthm:2}}

We start from the definition of $\hgamma_2$ in \eqref{est:2}. From the basic inequality, we know 
\begin{align*}
\frac{1}{2}\|\hgamma_2\|^2 - \LD \hgamma_2, \hgamma_1\RD + \lambda\|\hgamma_2\|_1\leq \frac{1}{2}\|\ttgamma\|^2 - \LD \ttgamma, \hgamma_1\RD + \lambda\|\ttgamma\|_1.
\end{align*}
Define $\bDelta = \hgamma_2 - \ttgamma$ and $\omega = \text{supp}(\ttgamma) = \text{supp}(\tgamma)$, then
\begin{align*}
\frac{1}{2}\|\hgamma_2 - \ttgamma\|^2\leq &\LD \bDelta,   \hgamma_1 - \ttgamma\RD + \lambda(\|\ttgamma\|_1 - \|\hgamma_2\|_1)\\
\leq &  \|\bDelta\|_1\|\ttgamma - \hgamma_1\|_{\infty} + \lambda\|\ttgamma_{\omega}\|_1 - \lambda\|(\hgamma_2)_\omega\|_1 - \lambda \|\bDelta_{\omega^c}\|_1\\
\leq &\|\bDelta\|_1\|\ttgamma - \hgamma_1\|_{\infty} + \lambda\|\bDelta_{\omega}\|_1 - \lambda \|\bDelta_{\omega^c}\|_1.
\end{align*}
By \eqref{b:15} and the setup of $\lambda$, we know $\|\ttgamma - \hgamma_1\|_{\infty}\leq \lambda/2$ with probability $1-\delta$. Therefore,
\begin{equation}\label{b:18}
\begin{aligned}
&\frac{1}{2}\|\bDelta\|^2_2\leq \frac{3\lambda}{2}\|\bDelta_{\omega}\|_1 - \frac{\lambda}{2}\|\bDelta_{\omega^c}\|_1 \Longrightarrow \|\bDelta\|_2\leq 3\sqrt{s}\lambda, \\
& \|\bDelta\|_1 \leq 4\|\bDelta_{\omega}\|_1\leq 4\sqrt{s}\|\bDelta\|_2\leq 12s\lambda.
\end{aligned}
\end{equation}
For bounding $\dist(\hgamma_2, \tgamma)$, we follow the same derivation as in \eqref{b:17}. This completes the proof.

\subsection{Proof of Theorem \ref{Sthm:4}}

We apply the one-dimensional $\sin(\theta)$ theorem in Lemma \ref{aux2:lem:2}. Recall that  $\tU = 2\mu_2\tgamma{\tgamma}^T$ in \eqref{equ:3}. We have
\begin{align}\label{b:19}
\|\hU& - \tU\|_{\infty, \infty} \nonumber\\
\leq& \|\hU - \mEbf[(y-\hf(\xb^T\hbeta))^2H(\xb)]\|_{\infty,\infty} + \|\mEbf[(y-\hf(\xb^T\hbeta))^2H(\xb)] - \tU\|_{\infty, \infty}.
\end{align}
For the second term, by the independence of $\epsilon$ and $\xb$,
\begin{align*}
\|\mEbf[(y-\hf(\xb^T\hbeta))^2H(\xb)] - \tU\|_{\infty, \infty} = & \|\mEbf[(f(\xb^T\tbeta) - \hf(\xb^T\hbeta))^2H(\xb)]\|_{\infty, \infty}.
\end{align*}
By argument in \eqref{cond:7}, we have
\begin{align}\label{b:20}
P\bigg(\bigg\|\mEbf[(y-\hf(\xb^T\hbeta))^2H(\xb)] - \tU\bigg\|_{\infty, \infty}\geq \frac{\sqrt{M}\cdot\tte_{f, \delta}(n, 1)}{\delta}\bigg)\leq 2\delta.
\end{align}
For the first term, we proceed as \eqref{b:14} and apply Lemma \ref{aux1:lem:2}. For any $\delta>0$, there exist constants $N_{\delta}, \Upsilon_1$ such that if $n\geq N_{\delta}$ and $\tau = \Upsilon_1(\frac{nM}{\log(2d^2/\delta)})^{\frac{1}{6}}$, we have
\begin{align}\label{b:21}
P\bigg(\biggr\|\hU - \mEbf[(y-\hf(\xb^T\hbeta))^2H(\xb)]\biggl\|_{\infty,\infty}>28\sqrt{\frac{\Upsilon_1M\log(2d^2/\delta)}{n}}\bigg)<4\delta.
\end{align}
Combining \eqref{b:19}, \eqref{b:20} and \eqref{b:21} together and adjusting $\delta$ properly, we know, for some constant $\Upsilon$,
\begin{align}\label{b:22}
P\bigg(\|\hU - \tU\|_{\infty, \infty}\leq \Upsilon\big(\sqrt{\frac{M\log(d/\delta)}{n}} + \frac{\sqrt{M}\cdot\tte_{f, \delta}(n,1)}{\delta}\big)\bigg)\geq 1-\delta.
\end{align}
Since $\|\hU - \tU\|_2\leq d\|\hU - \tU\|_{\infty, \infty}$, we further have
\begin{align*}
\|\hU - \tU\|_{2}\leq \Upsilon\big(d\sqrt{\frac{M\log(d/\delta)}{n}} + \frac{\sqrt{M}d\cdot\tte_{f, \delta}(n, 1)}{\delta}\big)
\end{align*}
with probability at least $1-\delta$. Without loss of generality, we assume $\mu_2>0$. Otherwise we can replace $\tU$ by $-\tU$ and $\hU$ by $-\hU$, but the estimator in \eqref{est:4} does not change, since we extract the eigenvector of $\hU$ corresponding to the eigenvalue with the largest magnitude. To apply Lemma \ref{aux2:lem:2}, we need the leading eigenvalue of $\hU$ to be positive. This is guaranteed by Weyl's inequality. In fact, since $\lambda_{\max}(\tU) = 2\mu_2>0$, 
for $n$ large enough, Lemma \ref{aux2:lem:3} implies
\begin{align*}
\lambda_{\max}(\hU)\geq 2\mu_2 - \|\hU - \tU\|_2>0.
\end{align*}
Thus, by Lemma \ref{aux2:lem:2}, we finally obtain 
\begin{align*}
\min_{\iota = \pm 1}\|\iota \hgamma_3 - \tgamma\|_2\leq \frac{\sqrt{2}\Upsilon}{\mu_2}\bigg( d\sqrt{\frac{M\log(d/\delta)}{n}} + \frac{\sqrt{M}d\cdot\tte_{f, \delta}(n, 1)}{\delta}\bigg),
\end{align*}
with probability at least $1-\delta$, which completes the proof.

\subsection{Proof of Theorem \ref{Sthm:5}}

Let $\tV = \tgamma{\tgamma}^T$ and $\tU = 2\mu_2\tV$. Since $\tV$ is feasible for the optimization program \eqref{equ:4}, from the basic inequality we have
\begin{align*}
\LD \hV, \hU\RD -\lambda\|\hV\|_{1,1}\geq \LD \tV,\hU\RD - \lambda\|\tV\|_{1,1}.
\end{align*}
This is equivalent to
\begin{align*}
\LD \hV - \tV, \hU - \tU\RD - \lambda\|\hV\|_{1,1} + \lambda\|\tV\|_{1,1}\geq \LD \tU, \tV - \hV\RD.
\end{align*}
For the right hand side term, we can assume $\mu_2>0$ without loss of generality. Otherwise we replace $\tU$ by $-\tU$. We apply Lemma \ref{aux2:lem:4} and get
\begin{align*}
\LD \tU, \tV - \hV\RD\geq \mu_2\|\tV - \hV\|_F^2.
\end{align*}
Applying the elementwise Holder's inequality, we get
\begin{align}\label{b:23}
\|\hV - \tV\|_{1,1}\|\hU - \tU\|_{\infty,\infty} - \lambda\|\hV\|_{1,1}+\lambda\|\tV\|_{1,1} \geq \mu_2\|\tV - \hV\|_F^2.
\end{align}
Denote $\bDelta = \hV - \tV$ and $\omega = \text{supp}(\tgamma)\times \text{supp}(\tgamma)$. By \eqref{b:22} and setup of $\lambda$, we have $\|\hU - \tU\|_{\infty,\infty}\leq \lambda$ with probability at least $1-\delta$. Therefore, the left hand side of \eqref{b:23} can be further upper bounded as
\begin{equation*}
\lambda(\|\bDelta\|_{1,1} - \|\hV\|_{1,1} + \|\tV\|_{1,1})
= 
\lambda(\|\bDelta_{\omega}\|_{1,1} - \|\hV_\omega\|_{1,1} + \|\tV_\omega\|_{1,1})
\leq 2\lambda\|\bDelta_{\omega}\|_{1,1}\leq 2s\lambda\|\bDelta\|_F.
\end{equation*}
Plugging back into \eqref{b:23}, we have
\begin{align*}
\|\bDelta\|_2\leq \|\bDelta\|_F\leq \frac{2s\lambda}{\mu_2}.
\end{align*}
Using Lemma \ref{aux2:lem:2} and \ref{aux2:lem:3}, we obtain
\begin{align*}
\min_{\iota = \pm 1}\|\iota \hgamma_4 - \tgamma\|_2\leq \frac{4\sqrt{2}s\lambda}{\mu_2},
\end{align*}
which finishes the proof.

\subsection{Proof of Theorem \ref{Sthm:3}}

We can apply same derivations as in Theorem \ref{Sthm:2}. We only need show $\lambda\geq 2\|\barw - \LD \hbeta, \hw\RD\cdot\hbeta - \ttgamma\|_{\infty}$ and the result then follows from \eqref{b:18}. Let $\bargamma = \barw - \LD \hbeta, \hw\RD\cdot\hbeta$. By \eqref{equ:2}, we have 
\begin{align}\label{b:24}
\|\bargamma - \ttgamma\|_{\infty}\leq &\|\barw - \mE[y^2S(\xb)]\|_{\infty} + \|\LD \hbeta, \hw\RD\cdot\hbeta - 2\eta_1\tbeta\|_{\infty} \nonumber\\
\leq & \|\barw - \mE[y^2S(\xb)]\|_{\infty} + 2|\eta_1|\cdot\|\hbeta - \tbeta\|_{\infty} + |\LD\hbeta,\hw\RD - 2\eta_1|.
\end{align}
For the first term, Lemma \ref{aux1:lem:2} suggests that, by setting $\tau = (\frac{nM}{\log(2d/\delta)})^{\frac{1}{6}}$,
\begin{align}\label{b:25}
P\bigg(\|\barw - \mE[y^2S(\xb)]\|_{\infty} \leq 7\sqrt{\frac{M\log(2d/\delta)}{n}}\bigg)\geq 1-\delta.
\end{align}
The second term has the rate from Theorem \ref{Supthm:1}. For the third term,
\begin{align}\label{b:26}
|\LD \hbeta,\hw\RD - 2\eta_1| = &|\LD \hbeta,\hw\RD - \LD \tbeta, \mE[y^2S(\xb)]\RD|  \nonumber\\
\leq & |\LD \hbeta, \hw\RD - \LD \tbeta, \hw\RD| + |\LD \tbeta, \hw\RD - \LD \tbeta, \mE[y^2S(\xb)]\RD| \nonumber\\
\leq  & \|\hbeta - \tbeta\|_{\infty}\|\hw\|_1 + \|\tbeta\|_1\|\hw - \mE[y^2S(\xb)]\|_{\infty} \nonumber\\
\leq & \|\hbeta - \tbeta\|_{\infty}\|\hw - \mE[y^2S(\xb)]\|_1 +\|\mE[y^2S(\xb)]\|_1\|\hbeta - \tbeta\|_{\infty}\nonumber\\
&+ \|\tbeta\|_1\|\hw - \mE[y^2S(\xb)]\|_{\infty}.
\end{align}
As proved in Theorem \ref{Sthm:2}, if $\kappa =  14\sqrt{\frac{M\log(2d/\delta)}{n}}\geq 2\|\barw - \mE[y^2S(\xb)]\|_{\infty}$, then
\begin{align*}
\|\hw - \mE[y^2S(\xb)]\|_1\leq &24s\kappa,\\
\|\hw - \mE[y^2S(\xb)]\|_{\infty}\leq & \|\hw - \barw\|_{\infty} + \|\barw - \mE[y^2S(\xb)]\|_{\infty}\leq \frac{3}{2}\kappa.
\end{align*}
The first inequality is from \eqref{b:18} noting that $\|\mE[y^2S(\xb)]\|_0\leq 2s$; the second inequality is due to the fact that $\|\bv - \phi_{\kappa}(\bv)\|_{\infty}\leq \kappa$, $\forall \bv\in\mR^d$. Plugging back into \eqref{b:26} and noting 
that $\|\mE[y^2S(\xb)]\|_1\leq 2|\eta_1|\cdot\|\tbeta\|_1 +|\mu_1|\cdot\|\tgamma\|_1 $,
\begin{align}\label{b:27}
|\LD \hbeta,\hw\RD - 2\eta_1|\leq \underbrace{24s\kappa\|\hbeta- \tbeta\|_{\infty}}_{\text{small order}} + (2|\eta_1|\cdot\|\tbeta\|_1 +|\mu_1|\cdot\|\tgamma\|_1)\|\hbeta - \tbeta\|_{\infty} + \frac{3}{2}\|\tbeta\|_1\kappa.
\end{align}
Combining \eqref{b:24}, \eqref{b:25}, \eqref{b:27}, Theorem \ref{Supthm:1} and using $\|\hbeta - \tbeta\|_{\infty}\asymp \kappa$, we have
\begin{align*}
\|\bargamma - \ttgamma\|_{\infty}\lesssim \underbrace{\big((|\eta_1|+1)\cdot\|\tbeta\|_1 +|\mu_1|\cdot\|\tgamma\|_1\big)}_{C_{(\eta_1, \mu_1, \|\tbeta\|_1, \|\tgamma\|_1)}}\kappa.
\end{align*}
Therefore, $\lambda \geq 2\|\bargamma - \ttgamma\|_{\infty}$ and the proof now follows as in Theorem \ref{Sthm:2}.

\subsection{Proof of Theorem \ref{Sthm:6}}

Based on the proof of Theorem \ref{Sthm:5} (see \eqref{b:23}), we only need show the setup of $\lambda$ satisfies $\lambda\geq \|\barU - \tU\|_{\infty, \infty}$, then the result holds by following the same derivations. By \eqref{equ:5} and the definition of $\barU$,
\begin{align*}
\|\barU - \tU\|_{\infty, \infty}\leq \|\ttU - \mE[y^2H(\xb)]\|_{\infty, \infty} + \|\hbeta^T\ttU\hbeta\cdot\hbeta\hbeta^T - 2\eta_2\tbeta\tbetat\|_{\infty, \infty}.
\end{align*}
By Lemma \ref{aux1:lem:2}, with $\tau = (\frac{nM}{\log(2d^2/\delta)})^{\frac{1}{6}}$, we have
\begin{align}\label{b:28}
P(\|\ttU - \mE[y^2H(\xb)]\|_{\infty, \infty} >7\sqrt{\frac{M\log(2d^2/\delta)}{n}})<\delta.
\end{align}
For the second term, we have
\begin{align*}
\|\hbeta^T\ttU\hbeta\cdot\hbeta\hbeta^T - 2\eta_2\tbeta\tbetat\|_{\infty, \infty}\leq & |\hbeta^T\ttU\hbeta - 2\eta_2| + 2|\eta_2|\cdot\|\hbeta\hbeta^T - \tbeta\tbetat\|_{\infty, \infty}\\
\leq & |\hbeta^T\ttU\hbeta - 2\eta_2| + 4|\eta_2|\cdot\|\hbeta - \tbeta\|_{\infty}.
\end{align*}
Furthermore,
\begin{align*}
|\hbeta^T\ttU\hbeta - 2\eta_2| = & |\hbeta^T\ttU\hbeta - \tbetat\mE[y^2H(\xb)]\tbeta| = |\LD \ttU, \hbeta\hbeta^T\RD - \LD \mE[y^2H(\xb)], \tbeta\tbetat\RD|\\
\leq & |\LD \ttU - \mE[y^2H(\xb)], \hbeta\hbeta^T\RD| + | \LD \mE[y^2H(\xb)], \hbeta\hbeta^T-\tbeta\tbetat\RD|\\
\leq &\|\ttU - \mE[y^2H(\xb)]\|_{\infty, \infty}\|\hbeta\hbeta^T\|_{1,1} + \|\mE[y^2H(\xb)]\|_{1,1} \|\hbeta\hbeta^T - \tbeta\tbetat\|_{\infty,\infty}\\
\leq & \|\ttU - \mE[y^2H(\xb)]\|_{\infty, \infty}\|\hbeta\|_{1}^2 + 4(|\eta_2|\cdot\|\tbeta\|_1^2 + |\mu_2|\cdot\|\tgamma\|_1^2)\|\hbeta - \tbeta\|_{\infty}\\
\leq & 2\underbrace{\|\ttU - \mE[y^2H(\xb)]\|_{\infty, \infty}\|\hbeta - \tbeta\|_{1}^2 }_{\text{small order}} + 2\|\tbeta\|_1^2\|\ttU - \mE[y^2H(\xb)]\|_{\infty, \infty}\\
&+4(|\eta_2|\cdot\|\tbeta\|_1^2 + |\mu_2|\cdot\|\tgamma\|_1^2)\|\hbeta - \tbeta\|_{\infty}.
\end{align*}
Combining the above display with \eqref{b:28} and Theorem \ref{Supthm:1}, for some constant $\Upsilon$, we have
\begin{align*}
P\bigg(|\hbeta^T\ttU\hbeta - 2\eta_2|\leq \Upsilon((|\eta_2|+1)\cdot\|\tbeta\|_1^2 + |\mu_2|\cdot\|\tgamma\|_1^2)\sqrt{\frac{M\log(2d^2/\delta)}{n}}\bigg)\geq 1-\delta.
\end{align*}
Hence, 
\begin{align*}
P\bigg(\|\barU - \tU\|_{\infty, \infty}\leq \underbrace{\Upsilon((|\eta_2|+1)\cdot\|\tbeta\|_1^2 + |\mu_2|\cdot\|\tgamma\|_1^2)}_{C'_{(\eta_2, \mu_2, \|\tbeta\|_1, \|\tgamma\|_1)}}\sqrt{\frac{M\log(2d^2/\delta)}{n}}\bigg)\geq 1-\delta.
\end{align*}
This concludes the proof.

\subsection{Proof of Theorem \ref{Mthm:1}}

From \eqref{Mest:1}, $\hG_1\in\mR^{d\times v}$ is a matrix whose columns are eigenvectors of $\hU$ corresponding to the top-$v$ eigenvalues. Based on \eqref{equ:6}, $\tU = 2\tG\bLambda\tGt$. Suppose $2\bLambda = P\Pi P^T$ to be the eigenvalue decomposition, then $\tU = \tG P\Pi P^T\tGt$. By \eqref{b:22}, we know, for any $\delta>0$, there exist $N_{\delta}$, $\Upsilon$ such that $\forall n\geq N_{\delta}$, 
\begin{equation}\label{b:29}
\begin{aligned}
\|\hU - \tU\|_2\leq & \|\hU - \tU\|_F\leq \Upsilon(d\sqrt{\frac{M\log(d/\delta)}{n}} + \frac{\sqrt{M}d\cdot\tte_{f, \delta}(n, 1)}{\delta})
\end{aligned}
\end{equation}
with probability at least $1-\delta$. Let $\lambda_v(\hU)$ be the $v$-th largest eigenvalue of $\hU$. According to Lemma \ref{aux2:lem:3}, when $n$ large enough,
\begin{align*}
\lambda_v(\hU)\geq &2\mu_3 - \|\hU - \tU\|_2>\frac{3\mu_3}{2},\\
\lambda_{v+1}(\hU)\leq & \|\hU - \tU\|_2<\frac{\mu_3}{2}.
\end{align*}
Therefore, the problem \eqref{Mest:1} has a unique solution up to an orthogonal transformation. Let $F = \tG P$ and note that $F^TF = I_v$, then
\begin{align*}
\LD \hU, FF^T\RD\leq \LD\hU, \hG_1\hG_1^T\RD.
\end{align*}
By Lemma \ref{aux1:lem:3} and \ref{aux2:lem:5}, we have
\begin{align*}
\inf_{Q\in\mQ^{v\times v}}\|\hG_1 - \tG Q\|_F = &\inf_{Q\in\mQ^{v\times v}}\|\hG_1 - \tG PQ\|\leq \sqrt{2}\|\sin(\angle(\hG_1, F))\|_F\\
\leq & \frac{1}{\mu_3}\|\hU - \tU\|_F\stackrel{\eqref{b:29}}{\leq} \frac{\Upsilon}{\mu_3}\bigg(d\sqrt{\frac{M\log(d/\delta)}{n}} + \frac{\sqrt{M}d\cdot\tte_{f, \delta}(n, 1)}{\delta}\bigg),
\end{align*}
which concludes the proof.

\subsection{Proof of Theorem \ref{Mthm:2}}

Define $\tV = \tG\tGt$ to be the projection matrix. Since $\tV$ is feasible for \eqref{equ:4}, we have
\begin{align*}
\LD \hV, \hU\RD -\lambda\|\hV\|_{1,1}\geq \LD \tV,\hU\RD - \lambda\|\tV\|_{1,1}.
\end{align*}
This implies
\begin{align*}
\LD \hV - \tV, \hU - \tU\RD - \lambda\|\hV\|_{1,1} + \lambda\|\tV\|_{1,1}\geq \LD \tU, \tV - \hV\RD.
\end{align*}
For the right hand side term, we apply Lemma \ref{aux2:lem:4} and get
\begin{align*}
\LD \tU, \tV - \hV\RD\geq \mu_3\|\tV - \hV\|_F^2.
\end{align*}
Same as \eqref{b:23} in the proof of Theorem \ref{Sthm:5}, as long as
\begin{align*}
\lambda\geq \Upsilon\bigg(\sqrt{\frac{M\log(d/\delta)}{n}} + \frac{\sqrt{M}\cdot\tte_{f, \delta}(n,1)}{\delta}\bigg)\stackrel{\eqref{b:22}}{\geq} \|\hU - \tU\|_{\infty, \infty},
\end{align*}
then
\begin{align*}
\|\hV - \tV\|_F\leq \frac{2s\sqrt{v}\lambda}{\mu_3}.
\end{align*}
Using Lemma \ref{aux2:lem:3}, we  see that $\hG_2$ is unique up to an orthogonal transformation and
\begin{align*}
\LD \hV, \tG\tGt\RD\leq \LD \hV, \hG_2\hG_2^T\RD.
\end{align*}
From Lemma \ref{aux1:lem:3} and \ref{aux2:lem:5}, we finally have
\begin{align*}
\inf_{Q\in\mQ^{v\times v}}\|\hG_2 - \tG Q\|_F\leq \sqrt{2}\|\sin(\angle(\hG_2, \tG))\|_F\leq 2\|\hV - \tV\|_F\leq \frac{4s\sqrt{v}\lambda}{\mu_3}.
\end{align*}
We finish the proof.

\subsection{Proof of Theorem \ref{Mthm:3}}

Same as the proof of Theorem \ref{Mthm:2}, as long as the setup of $\lambda$ satisfies $\|\barU - \tU\|_{\infty,\infty}\leq \lambda$, the result holds by following the same derivations. From the proof of Theorem \ref{Sthm:6},
\begin{align*}
\lambda\geq \Upsilon((|\eta_2|+1)\cdot\|\tbeta\|_1^2 + |\mu_2|\cdot\|\tgamma\|_1^2)\sqrt{\frac{M\log(2d^2/\delta)}{n}}
\end{align*}
satisfies the requirement. Thus, we conclude the proof.

\section{Auxiliary Lemmas}\label{sec:C}

We summarize lemmas used in Section \ref{sec:B}. Some of the proofs are provided in Section~\ref{sec:D}.

\begin{lemma}[Boundedness of $f(\xb^T\hbeta)$]\label{aux1:lem:1}
	
Suppose Assumption \ref{ass:1} holds with $p\geq 6$, $f$ is Lipschitz continuous and $\hbeta$ is a consistent estimator of $\tbeta$. Then
\begin{align*}
\mEb[|f(\xb^T\hbeta)|^6] \lesssim M_6
\end{align*}
for any $n\geq N_{M_6}$, where $N_{M_6}$ depends on $M_6$ only.
	
\end{lemma}

\begin{lemma}\label{aux1:lem:2}
	
Suppose $\{x_i, y_i, z_i, w_i\}_{i=1}^n$ is a sequence of $n$ i.i.d.~samples distributed as $x_i\sim x$, $y_i\sim y$, $z_i\sim z$, and $w_i\sim w$. Define $\wxi = \Psi_{\tau}(x_i) = x_i\cdot \pmb{1}_{\{|x_i|\leq \tau\}}$ to be the truncated samples (similar for $\wyi, \wzi, \wwi$). Suppose there exists a constant $M$ such that
\begin{align*}
\mE[|x|^{2(p+q+r+s)}] \vee \mE[|y|^{2(p+q+r+s)}] \vee \mE[|z|^{2(p+q+r+s)}] \vee \mE[|w|^{2(p+q+r+s)}] \leq M
\end{align*}
for some non-negative integers $p, q, r, s$. Then, for any $\delta>0$ and $\tau = (\frac{nM}{\log(2/\delta)})^{\frac{1}{2(p+q+r+s)}}$, we have
\begin{align*}
|\frac{1}{n}\sum_{i=1}^{n}\wxi^p\wyi^q\wzi^r\wwi^s - \mE[x^py^qz^rw^s]|\leq 7\sqrt{\frac{M\log(2/\delta)}{n}}
\end{align*}
with probability at least $1-\delta$.
	
\end{lemma}

\begin{definition}[Principal angles between two spaces]\label{aux2:def:3}
	
Let $A, B\in\mR^{d\times r}$ such that $A^TA = B^TB = I_r$ and $A^TB = U\Sigma V^T$ is  the singular value decomposition. 	Let $\angle(A,B) \in \mR^{r\times r}$ be the diagonal matrix with $\angle(A, B)_i = \arccos(\Sigma_i)$. We call $\angle(A,B)$ the $r$ principal angles between two subspaces $\IM(A)$ and $\IM(B)$.
	
\end{definition}

\begin{lemma}\label{aux1:lem:3}
	
Suppose $A, B\in\mR^{d\times r}$ have orthonormal columns. Then
\begin{align*}
\inf_{Q\in\mQ^{r\times r}}\|A-BQ\|_F^2\leq 2\|\sin(\angle(A,B))\|_F^2 = \|AA^T - BB^T\|_F^2,
\end{align*}
where the principal angles $\angle(A,B)$ are defined in Definition \ref{aux2:def:3}. 
	
\end{lemma}

\begin{lemma}[One-dimensional Davis-Kahan $\sin(\theta)$ theorem, Theorem 5.9 in \cite{Vershynin2012Introduction}]\label{aux2:lem:2} 
	
Suppose $A\in\mR^{d\times d}$ is a positive semidefinite matrix. Let $(\lambda_1, \bgamma_1),\ldots, (\lambda_d, \bgamma_d)$ denote the pairs of eigenvalues-eigenvectors of $A$ ordered such that $\lambda_1\geq \ldots\geq \lambda_d$. For any matrix $B\in\mR^{d\times d}$ such that the leading eigenvalue is positive, let $\bmu_1 \in \arg\max_{\|\bmu\|_2\leq 1}\bmu^TB\bmu$. Then
\begin{align*}
\min_{\delta = \pm 1}\|\delta\bmu_1 - \bgamma_1\|_2\leq \sqrt{2}\sin(\angle(\bmu_1,\bgamma_1))\leq \frac{2\sqrt{2}}{\lambda_1 - \lambda_2}\|A-B\|_{2},
\end{align*} 
where $\angle(\bmu_1, \bgamma_1) = \arccos(|\bmu_1^T\bgamma_1|)$.
	
\end{lemma}

\begin{lemma}[Weyl's inequality, \cite{Weyl1912Das}]\label{aux2:lem:3}
	
Suppose we have $A = B+C$ for symmetric matrices $A, B, C\in\mR^{d\times d}$, and their eigenvalues are denoted as $a_i, b_i, c_i$ in descending order. Then we have
\begin{align*}
b_i + c_d\leq a_i\leq b_i + c_1.
\end{align*} 
	
\end{lemma}

\begin{lemma}[Curvature, Lemma 3.1 in \cite{Vu2013Fantope}]\label{aux2:lem:4}
	
Let $A$ be a symmetric matrix and $E$ be the projection matrix that projects onto the subspace spanned by the eigenvectors of $A$ corresponding to its $d$ largest eigenvalues $\lambda_1\geq \lambda_2\geq \ldots \geq \lambda_d$. If $\delta_A = \lambda_d - \lambda_{d+1}>0$, then 
\begin{align*}
\frac{\delta_A}{2}\|E-F\|_F^2\leq \LD A, E-F\RD
\end{align*}
for all $F$ satisfying $0\preceq F\preceq I$ and $\TR(F) = d$.
	
\end{lemma}

\begin{lemma}[Variational $\sin(\theta)$ theorem, Corollary 4.1 in \cite{Vu2013Minimax}]\label{aux2:lem:5}
	
Let $A\in\mR^{p\times p}$ be a positive semidefinite matrix and suppose its eigenvalues $\lambda_1\geq\lambda_2\geq \ldots \geq \lambda_p$ satisfy $\delta_A = \lambda_d - \lambda_{d+1}>0$ for some $d<p$. Let $Q_1\in\mR^{p\times d}$ be the matrix whose columns are the eigenvectors of $A$ corresponding to the $d$ largest eigenvalues. We denote $E = Q_1Q_1^T$. Furthermore, suppose matrix $Q_2\in\mR^{p\times d}$ has orthonormal columns and let $F = Q_2Q_2^T$. Then for any symmetric matrix $B$, if it satisfies
\begin{align*}
\LD B, E\RD\leq \LD B, F\RD,
\end{align*}
we have
\begin{align*}
\|\sin(\angle(Q_1, Q_2))\|_F\leq\frac{\sqrt{2}}{\delta_A}\|A-B\|_F.
\end{align*}
Here $\sin(\angle(Q_1, Q_2))$ is defined in Definition \ref{aux2:def:3}.
	
\end{lemma}

\section{Proofs of Other Theorems and Lemmas}\label{sec:D}

This section presents the proofs for results in the appendix and Section \ref{sec:C}.

\subsection{Proof of Theorem \ref{Supthm:2}}

The proof for classical nonparametric regression has been established in Theorem 1 in \cite{Fan1993Local}. For completeness, we present a proof for the index model \eqref{mod:3}, which involves additional techniques. Our starting point is
\begin{align}\label{b:1}
\mEbx[|\hf(\xb^T\hbeta) - f(\xb^T\tbeta)|^4]\leq 8\mEbx[|\hf(\xb^T\hbeta) - f(\xb^T\hbeta)|^4] + 8|f(\xb^T\hbeta) - f(\xb^T\tbeta)|^4.
\end{align}
For the second term, by Assumption \ref{ass:3} (a) we have
\begin{align}\label{b:2}
|f(\xb^T\hbeta) - f(\xb^T\tbeta)|^4\leq L_1^4|\xb^T\tbeta - \xb^T\hbeta|^4.
\end{align}
To deal with the first term, we first introduce some additional notations. Let $t = \xb^T\hbeta$, $t_i = \xb^T_i\hbeta$, $\forall i\in[n]$, $Y = (y_1;\ldots;y_n;)$, $Y_0=(f(\xb_1^T\tbeta);\ldots;f(\xb_{n}^T\tbeta))$, $Y_1 = (f(t_1);\ldots;f(t_{n}))$, $\br = Y-Y_0 = (\ttg(\xb_1)\epsilon_1;\ldots;\ttg(\xb_{n})\epsilon_{n})$, and
\begin{align*}
X_t = \begin{pmatrix}
1  & t_1-t \\
1 & t_2-t \\
\vdots & \vdots \\
1 & t_{n}-t 
\end{pmatrix}, \text{\ \ } W_t = \begin{pmatrix}
K_h(t_1-t)\\
& K_h(t_2-t)\\
& & \ddots\\
& & & K_h(t_n-t)
\end{pmatrix},
\end{align*}
for $K_h(\cdot) = K(\cdot/h)/h$ where $K(\cdot)$ is Gaussian kernel. For $l=0,1,2\ldots$, we define
\begin{align*}
S_{nl} = &\sum_{j=1}^n(t_j - t)^lK_h(t_j-t), \text{\ \ } \Xi_{nl} =  \sum_{j=1}^n(t_j - t)^lK_h^2(t_j-t)\sqrt{\rbeta(t_j)},\\
\bar{a}_l = &\int |x|^lK(x) dx<\infty, \text{\ \ \ } \bar{b}_l = \int |x|^{l}K^2(x) dx<\infty,
\end{align*}
and also use $a_l, b_l$ to denote the integral without using absolute value. Note that taking expectation conditional on $\hbeta, \xb$ means evaluation point, $t$, is fixed and randomness comes from $\{\xb_i, \epsilon_i\}_{i\in[n]}$, which are independent from $\hbeta, \xb$. In what follows, we will drop off the subscript of matrix $X$ and $W$, and let $\mEt[\cdot] = \mEbx[\cdot]$. In addition, we define a $2$-by-$2$ matrix $H = \begin{pmatrix}
S_{n2} & -S_{n1}\\
-S_{n1} & S_{n0}
\end{pmatrix}$, a scalar $S_n = S_{n2}S_{n0} - S_{n1}^2 + \frac{1}{n^2h^2}$, and a vector $\bl = WXH\eb_1/S_n\in \mR^{n}$, where 
$\eb_1\in\mR^2$ is the first canonical basis of $\mR^2$.

With the above notations, the estimator proposed in \cite{Fan1993Local} can be written as $\hf(t) = \bl^TY$. We further have the following decomposition
\begin{align}\label{b:3}
\hf(t) - f(t) = &\eb_1^THX^TWY/S_n - f(t) \nonumber\\
= & \eb_1^THX^TW(Y-Y_0)/S_n + \eb_1^THX^TW(Y_0-Y_1)/S_n +\eb_1^THX^TWY_1/S_n - f(t) \nonumber\\
\eqqcolon & \I_1 + \I_2 + \I_3.
\end{align}
We proceed to bound $\mEt[\I_1^4]$, $\mEt[\I_2^4]$, $\mEt[\I_3^4]$. Since $\I_1 = \bl^T\br$, we have
\begin{align}\label{b:4}
\mEt[\I_1^4] =& \mEt[|\sum_{i\in[n]}\bl_i\br_i|^4] = \mEt[\sum_{i,j,s,r\in[n]} \bl_i\bl_j\bl_s\bl_r\ttg(\xb_i)\ttg(\xb_j)\ttg(\xb_s)\ttg(\xb_r)\epsilon_i\epsilon_j\epsilon_s\epsilon_r] \nonumber\\
= & \mE[\epsilon^4]\cdot\sum_{i\in[n]} \mEt[|\bl_i\ttg(\xb_i)|^4] + 6(\mE[\epsilon^2])^2 \sum_{i\neq j\in[n]}\mEt[\bl_i^2\bl_j^2\ttg(\xb_i)^2\ttg(\xb_j)^2]\nonumber\\
\leq & 6\mE[\epsilon^4] \big(\sum_{i\in[n]}\mEt[\bl_i^4\rbeta(t_i)] + \sum_{i\neq j\in[n]}\mEt[\bl_i^2\bl_j^2\sqrt{\rbeta(t_i)\rbeta(t_j)}]\big) \nonumber\\
= & 6\mE[\epsilon^4]\mEt[(\bl^T D_1 \bl)^2],
\end{align}
where $D_1 = \diag(\sqrt{\rbeta(t_1)}, \ldots, \sqrt{\rbeta(t_n)})$. For the term $\bl^TD_1\bl$, by simple calculations and we have
\begin{align}\label{b:5}
\bl^TD_1\bl = \frac{1}{S_n^2}(S_{n2}^2\Xi_{n0} - 2S_{n1}S_{n2}\Xi_{n1} + S_{n1}^2\Xi_{n2}).
\end{align}
For any random variable $Z_n$ and integer $r$, we write $Z_n = O_r(a_n)$, if $\mE[|Z_n|^r] = O(a_n^r)$, and define $o_r(a_n)$ similarly. Note that $Z_n = \mE[Z_n] + O_r\big((\mE[|Z_n - \mE[Z_n]|^r])^{1/r}\big)$ and, by Cauchy-Schwarz inequality, we have $O_r(a_n)O_r(b_n) = O_{r/2}(a_nb_n)$. We show how to control $S_{nl}$ in \eqref{b:5}, while the other terms are bounded in a similar way. Note that
\begin{align*}
\mEt[\frac{1}{nh^l}S_{nl}]= &\frac{1}{n}\sum_{j=1}^{n}\mEt[(\frac{t_j - t}{h})^lK_h(t_j-t)]= \int x^lK(x)q_{\hbeta}(t+hx)\, dx = a_l\cdot \qbeta(t) + O(h^{\alpha}).
\end{align*}
For any even integer $r\geq 2$, 
\begin{align*}
\mEt\bigg[\bigg|\frac{1}{nh^l}S_{nl} -& \mEt[\frac{1}{nh^l}S_{nl}]\bigg|^r\bigg] \\
= & \mEt\bigg[\bigg|\frac{1}{n}\sum_{j\in[n]}\bigg(\underbrace{(\frac{t_j - t}{h})^lK_h(t_j-t) - \mEt[(\frac{t_j - t}{h})^lK_h(t_j-t)]}_{\zeta_j}\bigg)\bigg|^r\bigg].
\end{align*}
For any positive integer $\ttr$, we know
\begin{align*}
\mEt[|\zeta_j|^\ttr] \leq 2^\ttr\mEt[|\frac{t_j - t}{h}|^{\ttr l}K_h^\ttr(t_j-t)] \leq  \frac{C_{\ttr, l}2^\ttr}{h^{\ttr-1}},
\end{align*}
where $C_{\ttr, l}$ is a constant depending on $\ttr, l$ and upper boundedness of $\qbeta$. Expanding the summation term and combining with above inequality, we can get that
\begin{align*}
\mEt[|\frac{1}{nh^l}S_{nl} - \mEt[\frac{1}{nh^l}S_{nl}]|^r] \leq&  \frac{1}{n^r} \mEt[(\sum_{j\in[n]}\zeta_j)^r]\\
= & \frac{1}{n^r} \sum_{k=1}^{r/2}\sum_{j_1,\ldots,j_k\in[n]}\sum_{\substack{c_1+\ldots+c_k = r\\c_i\geq 2, \forall i}}  \mEt\big[|\zeta_{j_1}|^{c_1}|\zeta_{j_2}|^{c_2}\cdots|\zeta_{j_k}|^{c_k}\big] \\
\lesssim & \frac{1}{n^r} \sum_{k=1}^{r/2}\begin{pmatrix}[]
n\\
k
\end{pmatrix}\begin{pmatrix}[]
r - k-1\\
k-1
\end{pmatrix}\frac{2^r}{h^{r-k}}\\
\leq &\max_{x\in\{2,3,\ldots,r/2\}}\begin{pmatrix}[]
r - x-1\\
x-1
\end{pmatrix}\frac{r2^r}{(nh)^{r/2}}.
\end{align*}
We can do similar calculation for $\Xi_{nl}$. Therefore, for any even integer $r>0$,
\begin{equation}\label{b:6}
\begin{aligned}
\frac{1}{nh^l}S_{nl} = &a_l\qbeta(t) + O_r(h^{\alpha} + \frac{1}{\sqrt{nh}}),\\
\frac{h}{nh^l}\Xi_{nl} = & b_l\qbeta(t)\sqrt{\rbeta(t)} + O_r(h^{\alpha}+ \frac{1}{\sqrt{nh}}). 
\end{aligned}
\end{equation}
Using \eqref{b:6}, we see the numerator in \eqref{b:5} can be written as
\begin{equation}\label{b:7}
\begin{aligned}
S_{n2}^2\Xi_{n0} - 2S_{n1}S_{n2}\Xi_{n1} + S_{n1}^2\Xi_{n2} = & n^3h^3a_2^2b_0\qbeta^3(t)\sqrt{\rbeta(t)}(1+ O_r(h^{\alpha}+ \frac{1}{\sqrt{nh}})).
\end{aligned}
\end{equation}
From (6.6) in \cite{Fan1993Local}, we also have
\begin{align}\label{b:8}
\frac{n^2h^2}{S_n} = \frac{1}{a_2\qbeta^2(t)} + o_4(1).
\end{align}
By \eqref{b:5}, \eqref{b:7}, \eqref{b:8}, there exists constant $\Upsilon_1$ such that
\begin{align*}
(\bl^TD_1\bl)^2 = \Upsilon_1\bigg(\frac{b_0^2\rbeta(t)}{n^2h^2\qbeta^2(t)} + \frac{o_1(1)}{n^2h^2}\bigg)\bigg(1+O_r(h^{\alpha} + \frac{1}{\sqrt{nh}})\bigg).
\end{align*}
Combining with \eqref{b:4} and let $n$ sufficient large to ignore the smaller order term, we have
\begin{align}\label{b:9}
\mEt[\I_1^4] \lesssim \mEt[(\bl^TD_1\bl)^2]\lesssim \frac{\rbeta(t)}{n^2h^2\qbeta^2(t)}.
\end{align}
Here we use the condition $h\rightarrow 0$ and $nh\rightarrow \infty$. For the term $\I_2$, by Taylor expansion
\begin{align*}
Y_0 - Y_1 = & (f'(t_1)\xb_1^T(\tbeta - \hbeta); \ldots; f'(t_n)\xb_{n}^T(\tbeta - \hbeta)) 
\\ 
& \qquad + (\frac{f''(\xi_1)(\xb_1^T(\tbeta - \hbeta))^2}{2};\ldots;\frac{f''(\xi_n)(\xb_n^T(\tbeta - \hbeta))^2}{2})\\
\eqqcolon &\bb_1 + \bb_2,
\end{align*}
where random variable $\xi_{j}\in (\xb_j^T\tbeta, \xb_j^T\hbeta)$. Therefore,
\begin{align*}
\mEt[\I_2^4] = \mEt[(\bl^T\bb_1 + \bl^T\bb_2)^4]\leq 8\mEt[(\bl^T\bb_1\bb_1^T\bl)^2] + 8\mEt[(\bl^T\bb_2\bb_2^T\bl)^2].
\end{align*}
We bound the first term as an example, while the second term can be easily shown to have smaller order error using the boundedness of $f''$. We have
\begin{align*}
(\bl^T\bb_1)^4 
&= \big(\sum_{i\in[n]} \bl_i f'(t_i)\xb_i^T(\tbeta - \hbeta)\big)^4 \leq \big(\sum_{i\in[n]}\bl_i^2(f'(t_i))^2\big)^2\big(\sum_{i\in[n]}(\xb_i^T(\tbeta - \hbeta))^2\big)^2\\
&\eqqcolon n^2(\bl^TD_2\bl)^2 \cdot \bigg(\frac{1}{n}\sum_{i\in[n]}(\xb_i^T(\tbeta - \hbeta))^2\bigg)^2,
\end{align*}
where $D_2 = \diag(|f'(t_1)|^2,\ldots,|f'(t_n)|^2)$. Using the condition that $\xb$ has bounded $4$-th moment in any direction, we know
\begin{align*}
\frac{1}{n}\sum_{i\in[n]}(\xb_i^T(\tbeta - \hbeta))^2 = O_2(\|\hbeta - \tbeta\|_2^2).
\end{align*}
Moreover, using similar derivations as bounding $(\bl^TD_1\bl)^2$, we obtain
\begin{align*}
(\bl^T\bb_1)^4 = \Upsilon_2\bigg(\frac{b_0^2(f'(t))^4O_{1}(\|\hbeta - \tbeta\|_2^4)}{h^2\qbeta^2(t)} + \frac{o_1(\|\hbeta - \tbeta\|_2^4)}{h^2}\bigg)\bigg(1+O_{r}(h^\alpha + \frac{1}{\sqrt{nh}})\bigg).
\end{align*}
Using the convergence rate of $\hbeta$, which implies $\mEt[o_1(\|\hbeta - \tbeta\|_2^4)/h^2]\rightarrow 0$, and ignoring the smaller order terms, we have
\begin{align}\label{b:10}
\mEt[\I_2^4]\lesssim \mEt[(\bl^T\bb_1)^4] \lesssim \frac{\|\hbeta - \tbeta\|_2^4|f'(t)|^4}{h^2\qbeta^2(t)}.
\end{align}
Lastly, we deal with $\I_3$ term. A simple observation is that
\begin{align*}
Y_1 = f(t)X\eb_1+f'(t)X\eb_2 + \bu,
\end{align*}
where $\bu \in\mR^{n}$ with $\bu_j = f(t_j) - f(t) - f'(t)(t_j-t) = \frac{f''(\xi_j)}{2}(t_j-t)^2$. Based on this, we have
\begin{align*}
\I_3 &= \eb_1^THX^TWY_1/S_n - f(t)= -\frac{f(t)}{n^2h^2S_n} +  \eb^T_1HX^TW\bu/S_n\\
= &-\frac{f(t)}{n^2h^2S_n}+\frac{1}{2S_n}\big(S_{n2}\sum_{j\in[n]}f''(\xi_j)(t_j-t)^2K_h(t_j-t) - S_{n1}\sum_{j\in[n]}f''(\xi_j)(t_j-t)^3K_h(t_j-t)\big)\\
\leq &-\frac{f(t)}{n^2h^2S_n}+\frac{L_1}{2S_n}(S_{n2}^2 + S_{n1}\sum_{j\in[n]}|t_j-t|^3K_h(t_j-t)).
\end{align*}
By simple calculations, we can obtain that
\begin{align*}
\sum_{j\in[n]}|t_j-t|^3K_h(t_j-t) = n\bar{a}_3\qbeta(t)(1+O_r(h^\alpha+\frac{1}{\sqrt{nh}})).
\end{align*}
Together with \eqref{b:6} and \eqref{b:8}, we have
\begin{align}\label{b:11}
\mEt[\I_3^4] \lesssim h^8 + \frac{|f(t)|^4}{n^{16}h^{16}\qbeta^8(t)}.
\end{align}
Combining results in \eqref{b:1}, \eqref{b:2}, \eqref{b:3}, \eqref{b:9}, \eqref{b:10} and \eqref{b:11}, we know there exists a constant $\Upsilon$ such that
\begin{align*}
\mEt[|\hf(t) - &f(t)|^4]\\
\leq& \Upsilon\bigg(h^8 + \frac{\rbeta(t)}{n^2h^2\qbeta^2(t)} +  \frac{\|\hbeta - \tbeta\|_2^4\cdot|f'(t)|^4}{h^2\qbeta^2(t)}  + |\xb^T\tbeta - \xb^T\hbeta|^4 + \frac{|f(t)|^4}{n^{16}h^{16}\qbeta^8(t)}\bigg)
\end{align*}
for a sufficiently large $n$. This concludes the first part of the proof and also gives an explicit formula for $\bare_f(\hbeta^T\xb, n, 1)$.

For the second part of proof, we will use the result in Corollary \ref{cor:1}. According to Assumption \ref{ass:3} (d), we directly have that for another constant $\Upsilon'$,
\begin{align*}
\hate_f(\hbeta, n, 1) = \mEb[\bare(\hbeta^T\xb, n, 1)]\leq \Upsilon'\bigg(h^8 + \frac{1}{n^2h^2} + \frac{\|\hbeta - \tbeta\|_2^4}{h^2}\bigg). 
\end{align*}
Plugging in the rate in Theorem \ref{Supthm:1}, we get
\begin{align*}
P\bigg(\sqrt{\hate_f(\hbeta, n, 1)}\gtrsim\underbrace{\big(h^4 + \frac{s\log (d/\delta)}{nh}\big)}_{\tte_{f, \delta}(n, 1)}\bigg)\leq \delta, \text{\ \ \ } \forall 0<\delta<1.
\end{align*}
Therefore, setting the bandwidth $h$ as $h\asymp n^{-1/5}$ we get $\tte_{f, \delta}(n, 1)\asymp n^{-4/5}$. Here we assume $s\log(d/\delta)$ is constant and negligible, otherwise we can choose the optimal bandwidth to be $h\asymp (s\log(d/\delta)/n)^{1/5}$ and get\footnote{In both cases $\tte_{f, \delta}(n, 1) = o(\sqrt{s\log d/n})$, hence, the dominant term in main theorems is always the parametric rate.}
\begin{align*}
\tte_{f, \delta}(n, 1) = \bigg(\frac{s\log(d/\delta)}{n}\bigg)^{4/5}.
\end{align*}
This implies condition (\ref{cond:6}) holds and we apply Corollary \ref{cor:1} to complete the proof.

\subsection{Proof of Lemma \ref{aux1:lem:1}}

By triangular inequality,
\begin{align*}
|f(\xb^T\hbeta)|^6\leq 32|f(\xb^T\hbeta) - f(\xb^T\tbeta)|^6 + 32|f(\xb^T\tbeta)|^6.
\end{align*}
By Assumption \ref{ass:1} and Lipschitz continuity of $f$, for any fixed $\hbeta$, we have
\begin{align*}
\mEb[|f(\xb^T\hbeta) - f(\xb^T\tbeta)|^6] \lesssim\|\hbeta - \tbeta\|_2^6.
\end{align*}
Using the consistency of $\hbeta$ and note that $\mEb[|f(\xb^T\tbeta)|^6] = \mE[|f(\xb^T\tbeta)|^6]\leq M_6$, we see
\begin{align*}
\mEb[|f(\xb^T\hbeta)|^6]\lesssim \|\hbeta -\tbeta\|_2^6 + M_6\lesssim M_6,
\end{align*}
for sufficient large $n$ that depends only on $M_6$.

\subsection{Proof of Lemma \ref{aux1:lem:2}}

Without loss of generality, we assume $\max(p,q,r,s)\geq 1$. We apply Bernstein's inequality in Corollary 2.11 in \cite{Boucheron2013Concentration}. We have
\begin{align}\label{d:1}
|\frac{1}{n}\sum_{i = 1}^{n}\wxi^p\wyi^q\wzi^r\wwi^s -& \mE[x^py^qz^rw^s]| \nonumber \\
\leq &|\frac{1}{n}\sum_{i = 1}^{n}\wxi^p\wyi^q\wzi^r\wwi^s  - \mE[\wx^p\wy^q\wz^r\ww^s]| + |\mE[\wx^p\wy^q\wz^r\ww^s] - \mE[x^py^qz^rw^s]| \nonumber\\ 
\eqqcolon &|\I_1| + |\I_2|.
\end{align}
For the term $\I_2$, we have
\begin{align}\label{d:2}
|\I_2|\leq & \mE[|x^py^qz^rw^s|\pmb{1}_{\{|x|>\tau \text{\ or\ } |y|>\tau \text{\ or\ } |z|>\tau \text{\ or\ } |w|>\tau\}}] \nonumber\\
\leq &\sqrt{\mE[x^{2p}y^{2q}z^{2r}w^{2s}]}\cdot\sqrt{P(|x|>\tau) + P(|y|>\tau)+ P(|z|>\tau) + P(|w|>\tau)} \nonumber\\
\leq & \sqrt{\big(\mE[x^{2(p+q+r+s)}]\big)^{\frac{p}{p+q+r+s}}\big(\mE[y^{2(p+q+r+s)}]\big)^{\frac{q}{p+q+r+s}}}\cdot \nonumber\\
&\sqrt{\big(\mE[z^{2(p+q+r+s)}]\big)^{\frac{r}{p+q+r+s}}\big(\mE[w^{2(p+q+r+s)}]\big)^{\frac{s}{p+q+r+s}}}\cdot\frac{\sqrt{4M}}{\tau^{p+q+r+s}} \nonumber\\
\leq &\frac{2M}{\tau^{p+q+r+s}}.
\end{align}
The third inequality is due to the generalized H\"older's inequality. For the term $\I_1$, we have
\begin{align*}
& |\wxi^p\wyi^q\wzi^r\wwi^s|\leq \tau^{p+q+r+s},\\
& V_n = \sum_{i=1}^{n}\VAR(\wxi^p\wyi^q\wzi^r\wwi^s)\leq n\mE[\wx^{2p}\wy^{2q}\wz^{2r}\ww^{2s}]\leq nM.
\end{align*}
By Bernstein's inequality, $\forall t>0$,
\begin{align}\label{d:3}
P(|\I_1|>t)\leq 2\exp(-\frac{nt^2}{2M + 4t\tau^{p+q+r+s}}).
\end{align}
For any $\delta>0$, we let $\frac{1}{\tau^{p+q+r+s}} = \sqrt{\frac{\log(2/\delta)}{nM}}$ and $t = 5\sqrt{\frac{M\log(2/\delta)}{n}}$. Plugging in \eqref{d:3} and combining with \eqref{d:1} and \eqref{d:2}, we obtain
\begin{align*}
\bigg|\frac{1}{n}\sum_{i = 1}^{n}\wxi^p\wyi^q\wzi^r\wwi^s - \mE[x^py^qz^rw^s]\bigg|\leq 7\sqrt{\frac{M\log(2/\delta)}{n}}
\end{align*}
with probability at least $1-\delta$. This concludes our proof.

\subsection{Proof of Lemma \ref{aux1:lem:3}}

Suppose $A^TB = U\cos(\angle(A,B))V^T$, then
\begin{align*}
\inf_{Q\in\mQ^{r\times r}}\|A-BQ\|_F^2 = &\inf_{Q\in\mQ^{r\times r}} (\|A\|_F^2 + \|B\|_F^2 - 2\TR(A^TBQ))\\
= & 2(r - \sup_{Q\in\mQ^{r\times r}}\TR(A^TBQ))\\
= & 2(r - \sup_{Q\in\mQ^{r\times r}}\TR(U\cos(\angle(A,B))V^TQ))\\
= & 2(r - \sup_{Q\in\mQ^{r\times r}}\TR(\cos(\angle(A,B))V^TQU)).
\end{align*}
Minimum is attained for $Q = VU^T$, hence
\begin{align*}
\inf_{Q\in\mQ^{r\times r}}\|A-BQ\|_F^2 = & 2(r - \sum_{i=1}^{r}\cos(\angle(A,B)_i))\\
\leq & 2(r - \sum_{i=1}^{r}\cos^2(\angle(A,B)_i))\\
= & 2\sum_{i=1}^{r}\sin^2(\angle(A,B)_i) = 2\|\sin(\angle(A,B))\|_F^2.
\end{align*}
Furthermore
\begin{align*}
\|\sin(\angle(A,B))\|_F^2 = & r - \|\cos(\angle(A,B))\|_F^2\\
= & \frac{1}{2}(\|AA^T\|_F^2 + \|BB^T\|_F^2 - 2\|A^TB\|_F^2)\\
= &\frac{1}{2}(\|AA^T\|_F^2 + \|BB^T\|_F^2 - 2\LD AA^T, BB^T\RD)\\
= & \frac{1}{2}\|AA^T - BB^T\|_F^2.
\end{align*}
Combining the two equations together, we complete the proof.

\end{document}